\documentclass[amscd,amssymb,verbatim,11pt]{amsart}
\numberwithin{equation}{section}
\usepackage{epsfig}
\usepackage{multido}
\usepackage{color}
\usepackage{pstricks}
\usepackage{pst-all}
\usepackage{pst-math}
\usepackage{pst-func}
\usepackage{pst-3dplot}
\usepackage{graphicx}
\usepackage{amsmath}
\usepackage{a4,amssymb}

\newtheorem{thm}{Theorem}[section]
\newtheorem{cor}[thm]{Corollary}
\newtheorem{lem}[thm]{Lemma}

\theoremstyle{definition}

\newtheorem{rem}[thm]{Remark}

\newif\ifShowLabels
\ShowLabelstrue
\newdimen\theight
\def\TeXref#1{
     \leavevmode\vadjust{\setbox0=\hbox{{\tt
            \quad\quad  {\small  \bf #1}}}%
     \theight=\ht0
     \advance\theight  by  \dp0
     \advance\theight  by  \lineskip
     \kern -\theight \vbox  to
     \theight{\rightline{\rlap{\box0}}%
      \vss}%
      }}%

    {\begin{thm}\label{#1} \ifShowLabels \TeXref{#1} \fi}%
    {\end{thm}}

    {\begin{def}\label{#1} \ifShowLabels \TeXref{#1} \fi}%
    {\end{def}}

    {\begin{lem}\label{#1} \ifShowLabels \TeXref{#1} \fi}%
    {\end{lem}}

    {\begin{cor}\label{#1} \ifShowLabels \TeXref{#1} \fi}%
    {\end{cor}}

\newcommand{\eqRef}[1]%
     {\ifShowLabels \TeXref{#1} \fi
      \begin{equation}\label{#1} }

\ShowLabelsfalse

\newcommand{\vsp}{\vskip 1em}
\newcommand{\vspp}{\vskip 2em}
\newcommand{\NI}{\noindent}
\newcommand{\bea}{\begin{eqnarray}}
\newcommand{\eea}{\end{eqnarray}}
\newcommand{\IR}{I\!\!R}
\newcommand{\bas}{\begin{align*}}
\newcommand{\eas}{\end{align*}}
\newcommand{\ba}{\begin{align}}
\newcommand{\ea}{\end{align}}

\newcommand{\be}{\begin{equation}}
\newcommand{\ee}{\end{equation}}
\newcommand{\ben}{\begin{eqnarray*}}
\newcommand{\een}{\end{eqnarray*}}
\newcommand{\lam}{\lambda}
\newcommand{\Om}{\Omega}
\newcommand{\om}{\omega}
\newcommand{\tht}{\theta}
\newcommand{\p}{\partial}
\newcommand{\al}{\alpha}

\newcommand{\ve}{\varepsilon}
\newcommand{\dl}{\delta}
\newcommand{\D}{\Delta}

\newcommand{\G}{\Gamma}

\newcommand{\s}{\sigma}

\newcommand{\lamo}{\lambda_{\Om}}
\newcommand{\lamb}{\bar{\lambda}}
\newcommand{\rh}{\rho}

\newcommand{\R}{\mathbb{R}}

\newcommand{\alh}{\hat{\al}}
\newcommand{\bh}{\hat{\beta}}
\newcommand{\gh}{\hat{\gamma}}
\newcommand{\kph}{\hat{\kappa}}
\newcommand{\kp}{\kappa}

\newcommand{\etb}{\bar{\eta}}

\newcommand{\phih}{\hat{\phi}}

\newcommand{\XO}{ \mathcal{P}^+_{\Om_T} }
\newcommand{\YO}{ \mathcal{P}^-_{\Om_T} }

\newcommand{\XC}{ \overline{\mathcal{P}}^+_{\Om_T} }
\newcommand{\YC}{ \overline{\mathcal{P}}^-_{\Om_T} }

\newcommand{\ex}{\exists}

\newcommand{\Df}{\D_\infty}

\title[{\large $\infty$}-Laplacian]{On the viscosity solutions to a degenerate parabolic differential equation}

\author[T. Bhattacharya and L. Marazzi]{Tilak Bhattacharya and Leonardo Marazzi}

\begin{document}

\maketitle

\section{Introduction}
 
\NI In this work, we study some properties of the viscosity solutions to a degenerate parabolic equation involving the non-homogeneous infinity-Laplacian. This may be viewed as related to the work done in \cite{BL}, where we studied an eigenvalue problem for the infinity-Laplacian. 
\vsp
\NI To make matters more 
precise, let $\Om\subset \IR^n,\;n\ge 2$, be a bounded domain and $\p\Om$ be its boundary. Define $\Om_T=\Om\times (0,T),\;T>0,$ and its parabolic boundary 
$P_T=(\Om\times\{0\} ) \cup ( \p\Om\times [0,T) )$. Our goal is to study existence and properties of positive solutions $\phi(x,t)$ to
\eqRef{pb1}
\Df \phi=(\phi^3)_t=3\phi^2 \phi_t,\;\mbox{in $\Om_T,$ with $\phi(x,0)=f(x),\;x\in \Om$ and $\phi |_{\p\Om\times [0,T)}=g(x,t),$}
\ee
where $f(x)$ and $g(x,t)$ are continuous. The operator $\Df$ is the non-homogeneous infinity-Laplacian and is defined as
$$\Df \phi=\sum_{i,j=1}^n \frac{\p \phi}{\p x_i} \frac{\p \phi}{\p x_j}\frac{\p^2 \phi}{\p x_i x_j}.$$
The infinity-Laplacian is nonlinear and is a degenerate elliptic operator. As a result, solutions are to be understood in the viscosity sense. Questions involving the infinity-Laplacian have been attracting considerable attention recently. In particular, existence, uniqueness and local regularity have been topics of great interest. For greater motivation and context, we direct the reader to the works \cite{ACJ,BDM,C,CEG, LS, LSA, PSSW}. Moreover, degenerate parabolic equations, in general, are of great interest, see for instance, \cite{ED}.
\vsp
\NI An additional reason for studying (\ref{pb1}) is the connection that it has with the following problem for the infinity-Laplacian, namely,
\eqRef{pb2}
\Df u(x)+\lam u(x)^3=0,\;\;\mbox{in $\Om$ with $u(x)=h(x),\;x\in \p\Om,$}
\ee
where $-\infty<\lam<\infty$, is a parameter, and $h\in C(\overline{\Om}).$
It is shown in \cite{BL} that the there is a $\lamo>0$ such that the problem (\ref{pb2}) has a solution $u\in C(\overline{\Om})$ for $0\le \lam<\lamo$. We refer to $\lamo$ as the first eigenvalue of $\Df$ on 
$\Om$.  The case $\lam<0$ has been discussed in \cite{BMO2}. Moreover, when $\lam=\lamo$ and $h=0$, then there is a positive eigenfunction $u$ that solves (\ref{pb2}). If we write 
$\phi(x,t)=v(x) w(t)$ then (\ref{pb1}) leads to
\eqRef{pb21}
\Df v=kv^3,\;\;\mbox{in $\Om$, and}\;\; 3w^{\prime}(t)=kw(t),\;\;0<t<\infty,
\ee
where $k$ is a constant. If $v$ is a solution then $\phi(x,t)=v(x) e^{kt/3}$ solves the partial differential equation in (\ref{pb1}). Thus, the problem in (\ref{pb2}) has close connections to
(\ref{pb1}). 
\vsp
\NI We make another observation. If $\phi>0$ solves the differential equation in (\ref{pb1}) and $\eta=\log \phi$ then a simple calculation yields
\eqRef{pb3}
\Df \eta+|D\eta|^4=3\eta_t,\;\;\;\mbox{in $\Om_T$ }.
\ee
We will provide full justification of this computation in Section 2. A similar change of variables for positive solutions $v$ to (\ref{pb21}) (set $w=\log v$) leads to
\eqRef{pb31}
\Df w+|Dw|^4=k,\;\;\mbox{in $\Om$}.
\ee
We feel that these questions may be of some interest in the context of the infinity-Laplacian. As a matter of fact the related equation $\Df \eta=\eta_t$ has been studied in 
\cite{AS}, \cite{AJK} and \cite{CW}. In \cite{AS}, the authors study the existence, uniqueness and regularity by employing a regularizing scheme. They show that solutions are Lipschitz continuous. In \cite{AJK}, the authors study asymptotic behavior of solutions to the same equation. The work in \cite{CW} presents a characterization of sub-solutions. 
We direct the reader to these works for a more in-depth discussion. All these works provide us the motivation to study (\ref{pb1}).
\vsp
\NI Our goals in this work is to prove a comparison principle, existence of positive solutions and some results about the asymptotic behavior of positive solutions. We do not employ any regularizations and work in the setting of viscosity solutions. For a more detailed discussion about viscosity solutions, see \cite{CIL}. 
Our work shows that positive solutions are continuous but provides no results on any additional regularity. It would be quite interesting to know if these are locally Holder continuous or even smoother. In this connection, it would be also interesting to know if a Harnack inequality holds for positive solutions. 
\vsp 
\NI We now describe the general lay out of the paper. In Section 2, we introduce definitions, notations, state some previously proven results relevant to our work. We prove also some preliminary results including the change of variables relating (\ref{pb1}) and (\ref{pb21}) to (\ref{pb3}) and (\ref{pb31}). We point out that the equation in (\ref{pb3}) plays an important role in our discussions. In Section 3, we prove a maximum principle that holds for any solution $\phi$ to (\ref{pb1}) regardless of its sign. Theorem \ref{comp} contains a comparison principle valid for positive solutions to (\ref{pb1}). The section ends with a partial result about the strong minimum principle. In Section 4, we discuss some partial results about large time asymptotic behavior of $\phi$ and provide some results on unbounded domains. We do not address any existence issues in this context.
Section 5 contains a proof of the existence of a positive continuous solution $\phi$ to (\ref{pb1}), see Theorem \ref{Exist14}. This leads to the existence of solutions to (\ref{pb3}).
Our approach is an adaptation of the ideas in \cite{CI, CIL} and uses the fact that one can find a sub-solution and a super-solution of (\ref{pb1}) that attain initial and boundary data. The construction of such sub-solutions and super-solutions is carried out in Appendix I (Section 6). This have been done since the work is quite computational in nature and, we feel, is best done in a separate section. Appendix II (Section 7) contains some calculations that have been used in Section 5.

\section{Notations and preliminary results}     
\vsp
\NI In this section we provide notations, definitions, some preliminary results and recall some useful lemmas proven previously in other works.
\vsp
\NI In what follows, $\Om\subset \IR^n,\;n\ge 2,$ will be a bounded domain, unless otherwise stated, and $\p\Om$ its boundary. Let $\overline{\Om}$ denote the closure of $\Om$ and
$\Om^c=\IR^n\setminus \Om$ denote the complement of $\Om$. The letters $x,\;y,\;z$ will often denote points in $\IR^n$. We reserve the letter $o$ for the origin in $\IR^n$. There will be occasion to write
a point $x=(x_1,x_2,\cdots, x_n)$. Moreover, for each $i=1,\cdots,n,$ $\vec{e_i}$ will denote the unit vector along the positive $x_i$ axis. The letters $s,\; t$ will denote points in $\IR^+\cup\{0\}$. For $T>0$, we define
\eqRef{def1}
\Om_T=\Om\times (0,T)=\{(x,t)\in \IR^n\times \IR:\;x\in \Om,\;0<t<T\}.
\ee
By $P_T$ we will denote the parabolic boundary of $\Om_T$ and this is the set 
\eqRef{def2}
P_T=(\Om\times \{0\}) \cup (\p \Om\times [0,T) ).
\ee
For future use, we set $S_T$ to be the side boundary, i.e., $S_T=\p\Om\times [0,T)$ and $I_\Om=\overline{\Om}\times \{0\}$, implying $P_T=I_\Om\cup S_T$.
In this work, we use $\Om_T\cup P_T$ and $\overline{\Om}\times[0,T)$ interchangeably.
Next, we describe open cylinders as follows. Let $B_r(x)\subset \IR^n$ be the ball of radius $r$ that is centered at $x$. For $r>0$ and $\tau>0$, we define the following cylinders
\eqRef{def3}
D_{r,\tau}(x,t)=B_r(x)\times(t-\tau/2, t+\tau/2),\;\;\;\mbox{and}\;\;F_{r,\tau}(x,t)=B_r(x)\times (t, t+\tau).
\ee  
\vsp
\NI For a function $u:\Om_T\cup P_T:\rightarrow \IR$, we recall the definitions of semi-jets $\mathcal{P}^+_{\Om_T}u(y,s)$ and $\mathcal{P}^-_{\Om_T}u(y,s)$ where $(y,s)\in \Om_T$, see \cite{CI, CIL} for a detailed discussion.
Let $\mathcal{S}(n)$ denote the set of symmetric $n\times n$ matrices. Then $(a,p,X)\in \IR\times \IR^n\times \mathcal{S}(n)$ is in $\mathcal{P}^+_{\Om_T}u(y,s)$ if
\eqRef{def4}
u(x,t)\le u(y,s)+a(t-s)+\langle p, x-y\rangle+\frac{ \langle X (x-y), x-y\rangle }{2}+o(|t-s|+|x-y|^2)
\ee
as $(x,t)\rightarrow(y,s)$, where $(x,t)\in \Om_T$. We define $\mathcal{P}^-_{\Om_T}u=-\mathcal{P}^+_{\Om_T}(-u).$ Next, 
$\XC u(y,s)$ is the set of all $(a,p,X)\in \IR\times \IR^n\times \mathcal{S}(n)$ such that for each $m=1,2,\cdots$, $\ex\;((x_m,t_m), a_m,p_m,X_m)\in \Om_T\times \IR\times\IR^n\times \mathcal{S}(n)$ with
$(a_m,p_m,X_m)\in \XO u(x_m,t_m)$ and $(x_m,t_m)\rightarrow (y,s)$, $u(x_m,t_m)\rightarrow u(y,s)$ and $(a_m,p_m,X_m)\rightarrow (a,p,X)$.
The set $\YC u(y,s)$ may now be analogously defined.
\vsp
\NI We discuss the notion of a viscosity solution of (\ref{pb1}) and (\ref{pb3}). Set
\eqRef{def5}
\Pi(\phi)=\Df \phi-3\phi^2 \phi_t\;\;\mbox{and}\;\; \Gamma (\eta)=\Df \eta+|D\eta|^4-3\eta_t.
\ee
We say that a function $u$ is a sub-solution of $\Pi(\phi)=0$ in $\Om_T$, if $u$ is upper semi-continuous in $\Om_T$ and 
\eqRef{def6}
\langle X p, p\rangle-3 a u^2 (x,t)\ge 0,\;\;\forall(x,t)\in \Om_T\;\mbox{and}\;\forall(a,p,X) \in \XO u(x,t).
\ee
We will write this as $\Pi(u)\ge 0$ in $\Om_T$. Similarly, a function $v$ is a super-solution of $\Pi(\phi)=0$, in $\Om_T$, if $v$ is lower semi-continuous in $\Om_T$ and
\eqRef{def7}
\langle Y q, q\rangle-3 b v^2 (x,t)\le 0,\;\forall(x,t)\in \Om_T\;\mbox{and}\;\forall(b,q,Y) \in \YO v(x,t).
\ee
Formally, we write $\Pi(v)\le 0.$ A function $u$ is a solution to $\Pi(u)=0$, if $u$ is continuous in $\Om_T$ and is both a sub-solution and a super-solution in $\Om_T$. Analogous definitions hold for the operator $\Gamma$, for instance, $v$ is a sub-solution to $\G(\eta)=0$, in $\Om_T$, if $v$ is upper semi-continuous in $\Om_T$ and 
\eqRef{def0}
\langle X p, p\rangle+|p|^4-3 a \ge 0,\;\;\forall(x,t)\in \Om_T\;\mbox{and}\;\forall(a,p,X) \in \XO v(x,t).
\ee
In this case, we write $\G(v)\ge 0$.
\vsp
\NI We say that $u$ is a sub-solution of (\ref{pb1}) if $u$ is upper semi-continuous in $\overline{\Om}\times [0,T)$, 
\eqRef{def8}
\Pi(u)\ge 0,\;\;\mbox{in $\Om_T$,}\;\; u(x,t)\le g(x,t),\;\forall(x,t)\in S_T\;\mbox{and}\;\;u(x,0)\le f(x),\;\forall x\in I_\Om.
\ee 
Similarly, $u$ is a super-solution of (\ref{pb1}) if $u$ is lower semi-continuous in $\overline{\Om}\times [0, T)$,
\eqRef{def9}
\Pi(u)\le 0,\;\;\mbox{in $\Om_T$,}\;\; u(x,t)\ge g(x,t),\;\forall\;(x,t)\in S_T\;\mbox{and}\;\;u(x,0)\ge f(x),\;\forall\;x\in I_\Om.
\ee
If both (\ref{def8}) and (\ref{def9}) hold then $u$ solves (\ref{pb1}), that is, $u\in C(\overline{\Om}\times[0,T))$, $\Pi(u)=0$ in $\Om_T$ and $u$ attains the initial data and the side condition. 
\vsp
\NI We now show that studying positive solutions to $\Pi(\phi)=0$ in $\Om_T$ is the same as studying solutions to $\G(\eta)=0$ in $\Om_T$, see (\ref{pb1}) and (\ref{pb3}). 
This equivalence will be used repeatedly in the work.
Before stating this result, we make an elementary observation. 

\begin{lem}\label{log}
Let $-1/3 \le c\le 1/3$. Then
\ben
0 \le \log(1+c)-\left(c-\frac{c^2}{2} \right) \le c^3,\;\;c\ge0,\;\;\mbox{and}\;\;\; c^3\le \log(1+c)- \left(c-\frac{c^2}{2} \right)\le 0,\;\;c<0.
\een
Also, 
$$0\le e^c- \left(1+c+\frac{c^2}{2}\right)\le c^3,\;\;c\ge 0,\;\;\mbox{and}\;\;\;c^3\le e^c-\left(1+c+\frac{c^2}{2} \right)\le 0,\;\;c<0.$$
\end{lem}
\NI{\bf Proof:} Consider the functions $f(c)=\log(1+c)-c+(c^2/2)-c^3$ and $g(c)=\log(1+c)-c+(c^2/2)$. Differentiating,
$f^{\prime}(c)=-c^2(2+3c)/(1+c)<0$ and $g^{\prime}(c)=c^2/(1+c)>0$, where $-1/3<c<\infty$. Since $f(0)=g(0)=0$, the inequalities
stated in the lemma hold. One can prove analogously the inequalities for $e^c$. $\Box$
\vsp
\NI We now prove a result about change of variables.

\begin{lem}\label{pde}
Let $T>0$ and suppose that $\phi:\Om_T\rightarrow \IR^+$ and $\phi\ge m$, for some $m>0.$ Then\\
\NI (i) $\Pi(\phi)\ge 0$ in $\Om_T$ if and only if $\eta=\log \phi$ satisfies $\G (\eta)\ge 0$ in $\Om_T$, and analogously,\\
\NI (ii) $\Pi(\phi)\le 0$ in $\Om_T$ if and only if $\eta=\log \phi$ satisfies $\G (\eta)\le 0$ in $\Om_T$.
\end{lem}
\NI{\bf Proof:} We prove (i). The proof of (ii) is similar. We recall (\ref{def6}). It is clear that $\phi$ is upper semi-continuous in $\Om_T$ if and only if $\eta$ is upper semi-continuous in $\Om_T$. 
\vsp
\NI Suppose that $\eta$ is a sub-solution, that is, $\G (\eta)\ge 0$ in $\Om_T$, see (\ref{def5}), (\ref{def6}) and (\ref{def0}). We will show that $\Pi(\phi)\ge 0$ in $\Om_T$.
Let $(y,s)\in \Om_T$ and $(a,p,X)\in \XO \phi(y,s)$. By (\ref{def4}),
\eqRef{pde1}
 \phi(x,t)\le \phi(y,s)+a(t-s)+\langle p, x-y \rangle+\frac{1}{2} \langle X (x-y), x-y \rangle +o(|t-s|+|x-y|^2)
 \ee
as $(x,t)\rightarrow(y,s)$, where $(x,t)\in \Om_T$.
Our goal is to show that
\eqRef{pde2}
\langle X p, p\rangle-3 a \phi^2 (y,s)\ge 0.
\ee
To do this we will show that there is a triple $(b,q,Y)\in \IR\times\IR^n\times \mathcal{S}(n)$, related to $(a,p,X)$, such that $(b,q,Y)\in \XO \eta(y,s).$
\vsp
\NI In (\ref{pde1}), writing $\phi=e^\eta$ and dividing both sides by $e^{\eta(y,s)}$ to obtain 
\eqRef{pde3}
\frac{e^{\eta(x,t)} }{e^{\eta(y,s)} }\le 1+  \frac{ a(t-s)+\langle p , x-y\rangle+\frac{ 1}{2}\langle X (x-y), x-y\rangle  }{e^{\eta(y,s)} } +o(|t-s|+|x-y|^2)
\ee
as $(x,t)\rightarrow (y,s)$, $(x,t)\in \Om_T$. Set
$$c= \frac{ a(t-s) + \langle p, x-y\rangle }{e^{\eta(y,s)} }+ \frac{1}{2} \sum_{i,j=1}^n \frac {X_{ij} }{e^{\eta(y,s) }} (x-y)_i (x-y)_j  +o(|t-s|+|x-y|^2),$$
where $(x,t)\rightarrow (y,s),\;(x,t)\in \Om_T$. In (\ref{pde3}), apply logarithm to both sides, use the second order expansion in Lemma \ref{log} and employ Young's inequality for terms of the type $|t-s|^a|x-y|^b,\;a,b\ge 1$, for instance, write
$$|t-s|^a|x-y|^b\le \frac{2|t-s|^{3a/2}}{3}+\frac{|x-y|^{3b}}{3}.$$
Thus, as $(x,t)\rightarrow (y,s)$, $(x,t)\in \Om_T$,
\ben
\eta(x,t)-\eta(y,s) \le  \log(1+c) &\le&   \frac{ a(t-s) + \langle p, x-y\rangle }{e^{\eta(y,s)} } + \frac{1}{2} \sum_{i,j=1}^n \frac{  X_{ij} (x-y)_i(x-y)_j } {e^{\eta(y,s)} } \\
&-&\frac{1}{2} \frac{  \langle p, x-y \rangle ^2 }{ e^{2\eta(y,s)} } +o(|t-s|+|x-y|^2).
\een
Clearly, as $(x,t)\rightarrow (y,s),\;(x,t)\in \Om_T$,
\ben
\eta(x,t)-\eta(y,s)  \le \frac{ a(t-s) + \langle p, x-y\rangle }{e^{\eta(y,s)} } &+& \frac{1}{2} \sum_{i,j=1}^n \left(  \frac {X_{ij} }{e^{\eta(y,s) }} - \frac{ p_ip_j}{ e^{2\eta(y,s)} } \right)(x-y)_i (x-y)_j  \\
&& +o(|t-s|+|x-y|^2),
\een
We set $b=ae^{-\eta(y,s)},\;q=pe^{-\eta(y,s)}$ and $Y_{ij}=X_{ij}e^{-\eta(y,s)}-p_ip_j e^{-2\eta(y,s)}$, and note that $(b,q,Y)\in \XO \eta(y,s)$. Since $\G(\eta)\ge 0$, we use (\ref{def0})
\ben
0\le \langle Y q, q\rangle+|q|^4-3 b =\sum_{i,j=1}^n \left( X_{ij} e^{-\eta(y,s)} - p_i p_j e^{-2\eta(y,s)} \right) p_i p_j e^{-2\eta(y,s)}+|p|^4e^{-4\eta(y,s)}  - 3  ae^{- \eta(y,s)}.
\een
Simplifying, we obtain
$$\langle X p, p\rangle  - 3 a e^{2\eta(y,s)}=\langle X p, p\rangle  - 3 a \phi^2(y,s)\ge 0,\;\;\forall(y,s)\in \Om_T,\;\;\forall(a,p,X)\in \XO \phi(y,s).$$
Thus, (\ref{pde2}) holds and $\Pi(\phi)\ge 0.$ 
\vsp
\NI Suppose now that $\Pi(\phi)\ge 0$ in $\Om_T$. Let $(y,s)\in \Om_T$ and $(b,q,Y)\in \XO \eta(y,s)$. 
We show that 
\eqRef{pde4}
\langle Y q,q\rangle+|q|^4-3 b\ge 0.
\ee
By (\ref{def4}), as $(x,t)\rightarrow (y,s)$, where $(x,t)\in \Om_T$, we have
$$\eta(x,t)\le \eta(y,s)+b(t-s)+\langle q, x-y\rangle+\frac{1}{2} \langle Y(x-y),x-y\rangle+o(|t-s|+|x-y|^2),$$
Writing $e^{\eta}=\phi$, 
$$
\phi(x,t)\le \phi(y,s) \exp \left( b(t-s)+\langle q, x-y\rangle+\frac{1}{2} \langle Y(x-y),x-y\rangle+o(|t-s|+|x-y|^2) \right)
$$
as $(x,t)\rightarrow (y,s)$. Using the second order expansion in Lemma \ref{log} and arguing as before, 
\ben
\phi(x,t)&\le& \phi(y,s)\left( 1+ b(t-s)+\langle q,x-y\rangle + \frac{1}{2} \sum_{i,j=1}^n  \left( Y_{ij}+ q_iq_j \right)(x-y)_i(x-y)_j   \right) \\
&&+o(|t-s|+|x-y|^2)
\een
as $(x,t)\rightarrow (y,s)$. Set $a=\phi(y,s)b,\;p=\phi(y,s) q$ and $X_{ij}=\phi(y,s) (Y_{ij}+q_iq_j).$ With this choice, $(a,p,X)\in \XO \phi(y,s) .$ Thus,
\ben
\langle Y q,q\rangle+|q|^4-3 b&=& \sum_{i,j=1}^n \left (\frac{X_{ij} }{\phi(y,s)} - q_iq_j \right)q_i q_j + |q|^4 - 3 \frac{a}{\phi(y,s)}\\
&=& \frac{ \langle X q, q\rangle  - 3  a}{\phi(y,s)} = \frac{ \langle X p, p\rangle  - 3  a \phi^2(y,s)}{\phi^3(y,s)}\ge 0.
\een
Thus, (\ref{pde4}) holds and $\Gamma(\eta)\ge 0$. $\Box$

\begin{rem}\label{pdes} The proof in Lemma \ref{pde} can be easily adapted to show equivalence for super-solutions. Thus, 
$\phi>0$ is a solution, i.e., $\Pi(\phi)=0$, in $\Om_T$ if and only if $\G(\eta)=0$, in $\Om_T$, where $\eta=\log \phi$.  
$\Box$
\end{rem}
\vsp
\NI We now list some results that have been proven in previous works. We start with the following versions of the comparison principle for the $\Df$. See Lemmas 2.2 and 2.3 in \cite{BMO2}.
\begin{lem}\label{comp1}
Suppose that $f\in C(\Om)$ and $f>0$, or $f<0$ or $f\equiv 0$ in $\Om$. Let $u,\;v\in C(\overline{\Om})$ satisfy $\Df u\ge f(x)$, in $\Om$ and $\Df v\le f(x)$, in $\Om$. Then
$$\sup_{\Om}(u-v)=\sup_{\p\Om}(u-v).$$
\end{lem}
\begin{lem}\label{comp2}
Suppose that $f_i:\Om\rightarrow \IR,\;i=1,\;2.$ Let $u,\;v\in C(\overline{\Om})$ satisfy $\Df u\ge f_1(x)$, in $\Om$ and $\Df v\le f_2(x)$, in $\Om$. Suppose that
$f_1(x)>f_2(x)$, for every $x\in \Om$ then
$$\sup_{\Om}(u-v)=\sup_{\p\Om}(u-v).$$ 
\end{lem}
\NI The next two results have been proven in \cite{BL}. For a definition of the first eigenvalue $\lamo$, see (3.2) and Theorem 4.1 in \cite{BL}.

\begin{thm}\label{exst}
Let $\Om\subset \IR^n$ be a bounded domain and $a\in C(\Om)\cap L^{\infty}(\Om)$, with $\inf_{\Om}a> 0$. Suppose that $h\in C(\Om)\cap L^{\infty}(\Om)$ and $b\in C(\p\Om)$. 
Then there is a number $\lamo>0$ such that if $0\le \lam<\lamo$ then there is a function $u\in C(\overline{\Om})$ that solves the following Dirichlet problem, 
\ben
\Df u+\lam a(x) u^3=h(x),\;\mbox{in $\Om,$ and $u=b$ on $\p\Om$.}
\een
If $h(x)=0$ and $\inf_{\p\Om} b>0$ then the above has a unique positive solution $u$. 
\end{thm}
\NI See Theorem 3.10, Lemma 3.3 and Remark 3.5 in \cite{BL}. Also, see (3.2) and Theorem 4.1 in \cite{BL}. 

\begin{lem}\label{bnd} 
Let $\Om\subset \IR^n$ be a bounded domain and $0\le \lam<\lamo$. Let $\dl>0$ and $u\in C(\overline{\Om}),\;u>0,$ solve the following Dirichlet problem, 
\ben
\Df u+\lam u^3=0,\;\mbox{in $\Om,$ and $u=\dl$ on $\p\Om$.}
\een
Then
$$\sup_{\Om} u\ge \dl \left( \frac{\lam}{\lamo-\lam} \right)^{1/3}.$$
\end{lem}
\NI See Theorem 4.1 and Lemma 3.3 in \cite{BL}. Theorem 4.1 in \cite{BL} also implies
\begin{thm}\label{eignv} 
Let $\Om\subset \IR^n$ be a bounded domain. Then there is a $\lamo>0$, the first eigenvalue of $\Df$, and $u\in C(\overline{\Om}),\;u>0,$ a first eigenfunction, that solve the following Dirichlet problem, 
\ben
\Df u+\lamo u^3=0,\;\mbox{in $\Om,$ and $u=0$ on $\p\Om$.}
\een
\end{thm}

\vsp
\NI We now recall a result about the radial version of Theorem \ref{exst}, on a ball, that will be used in this work, see Theorem 6.1 in \cite{BL}. Let $\Om=B_R(o),\;R>0$. We set $\lam_B=\lam_{B_R(o)}$ and $r=|x|.$ Let $u(x)=u(r),\;r=|x|$, then for $\dl\ge 0$, we consider
\eqRef{rad0}
\Df v+\lam v^3=\left(\frac{dv}{dr}\right)^2\frac{d^2v}{dr^2}+\lam  v(r)^3=0,\;\;\forall x\in B_R(o),\;\;\mbox{and}\;\;u(R)=\dl.
\ee

\begin{thm}\label{rad1}
Let $0<\lam<\lam_B$ and $\dl\ge 0$. Let $u>0$ solve
$$
u(x)=u(r)=m-(3\lam)^{1/3}\int_0^r \left[ \int_0^t  u(s)^3\;ds\right]^{1/3}\;dt,
$$
where $u(o)=m>0$ is so chosen that $u(R)=\dl$. Here, $\sup_{B_R(o)} u=m$.\\
\NI (i) If $\dl>0$ in (\ref{rad0}), then $u\in C(\overline{B_R(o)})$ and $u$ is the unique solution to (\ref{rad0}).\\
\NI (ii) If $\lam=\lam_B$, in (\ref{rad0}), then $u(R)=0$ and $u>0$ is a radial first eigenfunction. \\
\NI (iii) Also, $u$ defined by
$$
\int^1_{u(r)/m} \frac{ds}{(1-s^4)^{1/4} }=\lam ^{1/4} r, \;\; \forall\; 0<\lam\le \lam_B,
$$
is a radial solution to (\ref{rad0}) and we have, $\lam ^{1/4}R=\int^1_{\dl/m}\;(1-s^4)^{-1/4}\;ds.$
\end{thm}

\begin{rem}\label{rad2}
From part (iii) of Theorem \ref{rad1} with $\dl\ge 0$,
$$\int^1_{\dl/m} \frac{ds}{ (1-s^4)^{1/4} }=\lam^{1/4} R.$$
Fix $R$ in what follows. Clearly, $\dl/m$ is a decreasing continuous function of $\lam$.
Fixing $\dl$ and considering $m=m(\lam)$, we obtain 
\eqRef{rad3}
\frac{dm}{d\lam}=\frac{Rm^2}{4\dl\lam^{3/4}} \left( 1- \left ( \frac{\dl}{m} \right)^4 \right)^{1/4}>0.
\ee
Thus, $m$ is an increasing $C^1$ function of $\lam$. By Lemma \ref{bnd}, $m(\lam)\rightarrow \infty$ as $\lam\uparrow\lam_B$. 
Next, we fix $m$ and consider $\dl$ as a function of $\lam$. Thus, $\dl$ is continuous and decreasing, and Lemma \ref{bnd} implies that
$$\mbox{ $\dl\rightarrow 0$ as $\lam$ increases to $\lam_B$. }\;\;\Box$$ 
\end{rem}
\vsp
\NI We consider the problem in (\ref{rad0}) with $\lam<0$. Existence and uniqueness of a solution follow from Theorem 5.5 and Corollary 7.5 in 
\cite{BMO2}. 
\begin{lem}\label{rad4}
Let $u>0$ be continuous and defined in $[0,R]$, and $\dl> 0$. Fix $\lam>0$; if we define $u(r)$ by
\eqRef{rad5}
\int_\dl^{u(r)} \frac{ds}{ (s^4-\dl^4)^{1/4} }=\lam^{1/4} r,\;0\le r\le R,
\ee
then $u(r)\ge \dl$, in $0\le r\le R$, and is increasing in $r$. Writing $u(x)=u(|x|)$, $u$ solves
$$\Df u-\lam u^3=0,\;\;\mbox{in $B_R(o)$ and $u(o)=\dl.$} $$
Moreover, we have 
$$ u(r)=\dl+(3\lam)^{1/3} \int_0^r  \left( \int_0^t u(s)^3\;ds \right)^{1/3} dt,\;\;0\le r\le R.$$
\end{lem}
\NI{\bf Proof:} Set $r=|x|$. It is clear from (\ref{rad5}) that $u(o)=\dl$, $u$ is continuous on $[0,R]$ and $C^{\infty}$, in $r>0$. Also, $u$ is an increasing function, and in $r>0,$
\ben
\frac{du}{dr} =\lam^{1/4} \left( u(r)^4-\dl^4\right)^{1/4} ,\;\;\mbox{and}\;\;\frac{d^2u}{dr^2}=\lam^{1/4} \left( \frac{u(r)}{(u(r)^4-\dl^4)^{1/4}} \right)^3\;\frac{du}{dr}.
\een
This leads to
\eqRef{rad7}
\left(\frac{du}{dr}\right)^2\frac{d^2u}{dr^2}=\lam u^3,\;\;r>0.
\ee
Suppose that $\psi(x)$ is $C^2$ in $B_R(o)$ and $u(x)-\psi(x)$ has a maximum at $o$. Then
\eqRef{rad6}
0\le u(x)-u(o)\le \langle D\psi(o), x\rangle+\frac{ \langle D^2\psi(o)x, x \rangle}{2}+ o(|x|^2),\;\;\mbox{as $|x|\rightarrow 0.$}
\ee
Since $u(o)=\dl$ and $u(x)\ge \dl$, it follows quite easily that $D\psi(o)=0$. Next, we estimate, using $s\ge \dl$,
\eqRef{rad8}
\int_\dl^{u(r)} \frac{ds}{(s^4-\dl^4)^{1/4}}=\int_\dl^{u(r)} \frac{ds}{(s-\dl)^{1/4}(s^3+\dl s^2+\dl^2 s+\dl^3)^{1/4} }\le \frac{4^{3/4}(u(r)-\dl)^{3/4}}{3\dl^{3/4}}.
\ee
Using (\ref{rad5}), we find that
\ben
u(r)-\dl\ge \frac{3^{4/3} \dl}{4}\lam^{1/3} r^{4/3},\;\;r\ge 0.
\een
Using this on the left side of (\ref{rad6}),  
$$\frac{3^{4/3} \dl}{4}\lam^{1/3} |x|^{4/3}\le \frac{ \langle D^2\psi(o)x, x \rangle}{2}+ o(|x|^2),\;\;\mbox{as $|x|\rightarrow 0.$}.$$
Dividing by $|x|^2$ and letting $|x|\rightarrow 0$, one sees that $D^2\psi(o)$ does not exist. Thus, $u(x)$ is a sub-solution in $B_R(o).$ Now suppose that $\zeta\in C^2$ in $B_R(o)$ and $u(x)-\zeta(x)$ has a minimum at $o$, i.e.,
\eqRef{rad9}
u(x)-u(o)\ge \langle D\zeta(o),x\rangle +o(|x|),\;\;x\rightarrow 0.
\ee
Since $u(r)$ is continuous at $r=0$, we can make $u(r)\le 2\dl$, by taking $r$ small. Using $s=2\dl$ in the cubic polynomial in (\ref{rad8}), for small $r$, 
\ben
\frac{4}{3 (15^{3/4}) } \left(\frac{u(r)-\dl}{\dl} \right)^{3/4}\le \int_\dl^{u(r)} \frac{ds}{(s^4-\dl^4)^{1/4}}=\lam^{1/4}r.
\een
From (\ref{rad9}),
$$\langle D\zeta(o),x\rangle +o(|x|)\le  15\dl \lam^{1/3}|x|^{4/3},\;\;\mbox{as $|x|\rightarrow 0$}.$$
Clearly, $D\zeta(o)=0$ and $\Df \zeta(o)-\lam u(o)^3\le 0.$ The conclusion of the lemma holds. To show the last part, integrate (\ref{rad7}) and use the continuity of $u$ at $o$.
$\Box$

\begin{rem}\label{radp} Let $u$ be as in Lemma \ref{rad4}. We discuss an estimate for $u(R)$ in terms of $R$. From (\ref{rad5}), 
\eqRef{rad11}
\int^{u(R)}_\dl \frac{ds}{ (s^4-\dl^4)^{1/4} }=\lam^{1/4} R.
\ee
Then $u(R)=\sup_{B_R(o)} u$, and, if we fix $R$, then $u(R)$ is a $C^1$ increasing function of $\lam$ and grows unboundedly as $\lam$ increases. Similarly, if we fix $\lam$ then $u(R)$ grows unboundedly as $R$ increases. Take $R>0$, large, so that $u(R)\ge 2\dl$. To get a better estimate for the latter, we write
\eqRef{rad12}
 \int_\dl^{u(R)} \frac{ds}{ (s^4-\dl^4)^{1/4} }=\int_\dl^{2\dl}\frac{ds}{ (s^4-\dl^4)^{1/4} }+\int_{2\dl}^{u(R)}\frac{ds}{ (s^4-\dl^4)^{1/4} }.
 \ee
Set
$$A=\int_\dl^{2\dl} \frac{ds}{ (s^4-\dl^4)^{1/4} }=\int_1^2\frac{ds}{ (s^4-1)^{1/4} }.$$ 
In the second integral on the right side of (\ref{rad12}), we use the inequality $s^4(1-2^{-4})\le s^4-\dl^4\le s^4,\;s\ge 2\dl.$ Applying (\ref{rad5}),
$$A+\log\left( \frac{u(R)}{2\dl} \right)\le \lam^{1/4} R\le  A+ \frac{2}{(2^4-1)^{1/4}}  \log\left( \frac{u(R)}{ 2\dl} \right),\;\;\;\mbox{for large $R$.} $$
Set $\s=2/(2^4-1)^{1/4}$. Then there is a constant $B>0$ such that 
\eqRef{rad13}
2\dl B^{1/\s} e^{(\lam^{1/4}R/\s)}\le u(R)\le 2B\dl e^{\lam^{1/4}R},\;\;\mbox{for large $R$.}
\ee
This estimate will prove useful on unbounded domains. $\Box$ 
\end{rem}

\NI Next, we verify that if $u(x)$ is as in Theorem \ref{rad1} or as in Lemma \ref{rad4} then $\phi(x,t)=u(x)e^{kt}$, for an appropriate $k$, solves
$\Df \phi=3\phi^2 \phi_t$ in the sense of viscosity.

\begin{lem}\label{sol1}
Let $\Om\subset \IR^n$, be a bounded domain. Suppose that, for some $\lam\in \IR$, $u\in C(\overline{\Om}),\;u\ge 0,$ solve  
$$\Df u+\lam u^3=0,\;\;\mbox{in $\Om$ and $u=f(x)$ on $\p\Om$,}$$
where $f\in C(\p\Om)$. Let $T>0$, then for $3k\ge(\le)-\lam$, the function
$\phi(x,t)=u(x) e^{kt}$ solves
$$\Df \phi\le (\ge)3\phi^2 \phi_t,\;\mbox{in $\Om_T$,}$$
where $\phi(x,t)=f(x) e^{kt}$, $\forall(x,t)\in \p\Om \times [0,T),$ and $\phi(x,0)=u(x),\;\forall x\in \Om$.
\vsp
\NI If $3k=-\lam$, then no sign restriction on $u$ is necessary, and $\phi(x,t)=u(x)e^{kt}$ solves $\Df\phi=3\phi^2\phi_t$ in $\Om_T$.
\end{lem}
\NI{\bf Proof:} We will prove this for the case $k\le -\lam/3$. The case $k\ge -\lam/3$ will follow in an analogous manner.
Assume that $k\le -\lam/3$ and $u\ge 0.$ Let $(y,s)\in \Om_T$ and $(a,p,X)\in \XO \phi(y,s)$, then
\eqRef{sol01}
\phi(x,t)\le \phi(y,s)+a(t-s)+\langle p, x-y\rangle+ \frac{ \langle X(x-y), x-y\rangle }{2}+o(|t-s|+|x-y|^2)
\ee
as $(x,t)\rightarrow (y,s),\;(x,t)\in \Om_T.$ Since $\phi(x,t)=u(x) e^{kt}$, taking $t=s$ in (\ref{sol01}) and dividing by $e^{ks}$, we get
$$u(x)\le u(y)+\langle pe^{-ks}, x-y\rangle+ \frac{ \langle Xe^{-ks}(x-y), x-y \rangle }{2}+o(|x-y|^2), $$
as $x\rightarrow y$. Since $u$ is a solution to $\Df u+\lam u^3=0$, we obtain
$e^{-3ks} \langle Xp,p\rangle+\lam u(y)^3\ge 0,$
leading to
\eqRef{sol2}
\langle Xp, p\rangle+\lam e^{3ks}u(y)^3=\langle Xp,p\rangle+\lam \phi^3(y,s)\ge 0.
\ee
Next, in (\ref{sol01}) set $x=y$ to obtain
$u(y)(e^{kt}-e^{ks})\le a(t-s)+o(|t-s|),\;\;\mbox{as $t\rightarrow s$.}$
It follows that $a=ke^{ks}u(y)\le-\lam e^{ks} u(y)/3$. Since $\phi\ge 0$, using the expression for $a$ in (\ref{sol2}), we see that
\eqRef{sol3}
\langle X p, p\rangle\ge 3 a \phi^2(y,s).
\ee
Thus, $\phi(x,t)$ is a sub-solution. In the even $3k=\lam$, we have that $a=-\lam u(y) e^{-\lam s/3}/3$. The inequality in (\ref{sol3}) holds regardless of the sign of $\phi$ and implies
that $\Df \phi\ge 3\phi^2 \phi_t$ in $\Om_T$. One can now argue in a similar fashion that $\Df \phi\le 3\phi^2 \phi_t$. The lemma holds.
$\Box$
\vsp
\begin{rem}\label{sol4} Let $u(x),\;\forall \bar{\Om}$, be as in the statement of Lemma \ref{sol1}. Suppose that $g:[0,\infty)\rightarrow \IR$ is a $C^1$ function. Arguing as done 
in the proof of Lemma \ref{sol1}, one can show that if $\phi(x,t)=u(x)g(t),\;\forall(x,t)\in \Om_T$ then
$$\Df \phi-3\phi^2\phi_t=-g(t)^2u(x)^3\left( \lam +3g^{\prime}(t) \right),\;\;\;\forall(x,t)\in \Om_T.\;\;\;\Box$$
\end{rem}
\vsp

\section{ Maximum and comparison principles}
\vsp
\NI Our main goal in this section is to prove a maximum principle and a comparison principle for (\ref{pb1}). We have included both here since our proof of the comparison principle applies only to solutions that are positive or have one sign in $\Om_T$. It is not clear to us if a more general version holds.
The maximum principle, on the other hand, applies to solutions that may change sign.  
We also provide a partial result regarding the strong minimum principle.  
\vsp
\NI We start by proving the weak maximum principle for sub-solutions to (\ref{pb1}). We note that $\phi$ is a sub-solution if and only if $-\phi$ is a super-solution. This will imply the weak minimum principle. We introduce the following notation. For any $x\in \IR^n$, write $x=(x_1,x_2,\cdots, x_n)$ and define the projections $\pi_i(x)=x_i,\;i=1,2,\cdots,n.$ 

\begin{lem}\label{max}{(Weak Maximum Principle)} Suppose that $\phi$ is upper semi-continuous in $\Om_T\cup P_T$ and solves 
$$\Df \phi\ge 3\phi^2 \phi_t,\;\;\mbox{in $\Om_T$}.$$
Then $\sup_{\Om_{T}}\phi\le\sup_{P_T} \phi=\sup_{\Om_T\cup P_T}\phi.$ In the event, $\phi$ is lower semi-continuous in $\Om_T\cup P_T$ and $\Df \phi\le 3\phi^2\phi_t$ in $\Om_T$, 
then $\inf_{\Om_T} \phi\ge \inf_{P_T}\phi=\inf_{\Om_T\cup P_T}\phi.$
\end{lem}
\NI{\bf Proof:} Let $0<\tau<T$ with $\tau$ close to $T$. We show that the weak maximum principle holds in $\Om_\tau$ for any $\tau<T$. Fix $\tau<T$. Note that $\overline{\Om}_\tau$ is compactly contained in $\Om_T\cup P_T$.
\vsp
\NI We record the following for later use.
Since $\Om$ is bounded, there are $\rho,\;\bar{\rho}\in \IR$ such that $-\infty<\rho\le \pi_1(x)=x_1\le \bar{\rho}<\infty,\;\forall x\in \Om$. Set 
Choose $\ell>3$ such that $\ell R>\bar{\rho}-\rho$. Thus, 
$$0\le \pi_1(x)-\rho\le 2\ell R,\;\forall x\in \Om.$$
\NI Set 
\eqRef{max0}
M=\sup_{\Om_\tau} \phi,\;\;L=\sup_{P_\tau}\phi \;\;\mbox{and}\;\;C=\sup_{\overline{\Om}}\phi(x,\tau).
\ee
Call $\dl=M-L$ and assume that $\dl>0$. Now call
\eqRef{max1}
A=\max(M-L,\;C-L).
\ee
Since $\sup_{\Om_\tau} \phi >\sup_{P_\tau}\phi(x,t)$ and $\phi(x,t)$ is upper semi-continuous in $\overline{\Om}_\tau$, either there is a point $(y,s)\in \Om_\tau$ such that $\phi(y,s)=M$, or $\phi(x,t)<M,\;\forall(x,t)\in \Om_\tau$ and there is sequence $(y_k,s_k)\in \Om_\tau$ such that $s_k\uparrow \tau$ and $\phi(y_k,s_k)\rightarrow M$, as $k\rightarrow \infty$. We take up the second possibility first. Since $s_k\uparrow \tau$, we choose $0<s_k<\tau$ such that $\phi(y_k,s_k)>L+3\dl/4$ and $0<s_k<\tau$. Fix $s_k$
and define
$$g(t)= \left\{ \begin{array}{lcr} 0, && 0\le t\le s_k\\ \left( \frac{t-s_k}{\tau-s_k} \right)^4, && s_k\le t\le \tau. \end{array} \right. $$
Then $g\ge 0$, $g$ is $C^2$ and $g^{\prime}(t)\ge 0.$ Next, select $0<\ve\le \min(\dl/4,1/2),$ recall (\ref{max0}), (\ref{max1}) and set
\ben
\psi(x,t)=L+\frac{\ve}{4} +Ag(t)-\frac{\ve}{16} \left(1-\left(\frac{\pi_1(x)-\rho}{3\ell R} \right)\right)^2,\;\;\forall(x,t)\in \overline{\Om}_{\tau}.
\een
Then 
$\psi(x,t)\ge L+\ve/8,\;\forall(x,t)\in \overline{\Om}_\tau$ and $\psi(x,\tau)\ge C+\ve/8,\;\;\forall x\in \overline{\Om}$, and
\ben
\phi(y_k,s_k)-\psi(y_k, s_k)\ge L+\frac{3\dl}{4}-L-\frac{\ve}{4}=\frac{3\dl}{4}-\frac{\ve}{4}
\ge \frac{3\dl}{4}-\frac{\dl}{16}>\frac{\dl}{4}>0.
\een
Since $\phi-\psi\le 0$ on $\p \overline{\Om}_{\tau}$, $\phi-\psi$ has a positive maximum at some point $(z,\tht)\in \Om_{\tau}$ with $\tht>0$. Calculating and using (\ref{max0}), we get
$$\Df \psi=(D_1\psi)^2D_{11}\psi=-\frac{\ve^3}{8^3 3^4 \ell^4R^4}\left( 1- \left( \frac{x_1-\rho}{3\ell R} \right) \right)^2\le -\frac{\ve^3}{8^3 3^4 \ell^4R^4}\left(\frac{8}{9}\right)^2<0.$$
Thus,
$$\Df \psi(z,\tht)<0\le 3\phi(z,\tht)^2\psi_t(z,\tht).$$
We obtain a contradiction. For the possibility $\phi(y,s)=M$, where $(y,s)\in \Om_\tau$, we apply the above argument with $s_k=s$. 
\vsp
\NI The above implies that $\sup_{\Om_\tau}\phi\le\sup_{P_\tau}\phi$ for any $\tau<T.$ If $\sup_{\Om_T}\phi>\sup_{P_T} \phi$ then there is a point $(y,s)\in \Om_T$ (with $s<T$) such that 
$\phi(y,s)>\sup_{P_T} \phi$. Select $s<\bar{s}<T$. Then, $\sup_{P_T}\phi<\phi(y,s)\le \sup_{\Om_{\bar{s}}}\phi\le \sup_{P_{\bar{s}}}\phi\le \sup_{P_T}\phi.$ This is a contradiction and the lemma holds.
$\Box$
\begin{rem}\label{min} 
\NI Lemmas \ref{max} and \ref{pde} implies that the maximum principle holds for sub-solutions of $\Df \eta+|\eta|^4\ge 3\eta_t$. 
\vsp
\NI If $\phi\in C(\Om_T\cup P_T)$ then Lemma \ref{max} implies that $\sup_{\Om_T}\phi=\sup_{P_T}\phi.\;\;\;\Box$
\end{rem}
\vsp
\NI Our next goal in this section is to prove a comparison principle for (\ref{pb1}). Unlike the maximum principle proven in Lemma \ref{max}, we will require that the sub-solutions and super-solutions have one sign. There are simple situations where this is not needed. For instance, if $u$ is a sub-solution and $v$, a super-solution, and $\inf_{P_T} v\ge \sup_{P_T} u$ then by the maximum principle, $u\le v$ in $\Om_T$.

\begin{thm}\label{comp}{(Comparison principle)}
Suppose that $\Om\subset \IR^n$ is bounded and $T>0$. Let $u$ be upper semi-continuous and $v$ be lower semi-continuous in $\Om_T\cup P_T$, and satisfy
$$\Df u\ge 3 u^2 u_t,\;\;\mbox{and}\;\;\Df v\le 3 v^2 v_t,\;\;\mbox{in $\Om_T$}.$$
Assume that $\min(\inf_{\Om_T\cup P_T} u, \inf_{\Om_T\cup P_T} v)>0$. If $\sup_{P_T}v<\infty$ and
$u\le v$ on $P_T$, then $u\le v$ in $\Om_T$.
\end{thm}
\NI{\bf Proof:} We employ the ideas in \cite{CIL} after making the change of variables described in Lemma \ref{pde}, also see Remark \ref{pdes}. We note that since $u\le v$ on $P_T$, by Lemma \ref{max}, $u$ is bounded in $\Om_T$.
\vsp
\NI Define $\eta(x,t)=\log u(x,t)$ and $\zeta(x,t)=\log v(x,t).$ Then $\eta$ is bounded and $\zeta$ is bounded, in particular,  from below. By Lemma \ref{pde},
\eqRef{comp1}
\Df \eta+|D\eta|^4\ge 3 \eta_t\;\;\;\mbox{and}\;\;\;\Df \zeta+|D\zeta|^4\le 3\zeta_t\;\;\;\mbox{in $\Om_T$}
\ee
with $\eta\le \zeta$ on $P_T$.
Our goal is to prove the comparison principle for (\ref{comp1}). We record the following calculation for later use. For $\ve>0$ (to be determined later) 
if $w$ is a sub-solution of (\ref{comp1}), define $w_{\ve}(x,t)=w(x,t)-\ve/(T-t),\;\forall(x,t)\in \Om_T.$ Then
\eqRef{comp2}
\Df w_{\ve}+|Dw_{\ve}|^4\ge 3 (w_{\ve})_t + \frac{3 \ve}{(T-t)^2},\;\;\;\mbox{in $\Om_T$.}
\ee
\vsp
\NI Let $\dl=\sup_{\Om_T}(\eta-\zeta)$ and suppose that $\dl>0$. Then there is a point $(y,\tau)\in \Om_T$, such that $(\eta-\zeta)(y,\tau)>3\dl/4$. Choose $\ve$ such that $0<\ve<\dl(T-\tau)/8.$ Writing $\eta_\ve=\eta-(\ve/(T-t))$, we have
\eqRef{comp3}
\sup_{\Om_T} (\eta_\ve-\zeta)\ge \frac{\dl}{2}.
\ee
Fix $\ve>0$ and recall from (\ref{comp2}) that $\eta_{\ve}$ is a sub-solution of (\ref{comp1}) and 
\eqRef{comp4}
\eta_\ve(x,t)\le \eta(x,t)-\ve/T,\;\;\mbox{in $\Om_T,\;\;\eta_\ve\le \zeta-\ve/T$ on $P_T$, and}\;\mbox{$\lim_{t\uparrow T} (\sup_{\Om}\eta_\ve(x,t))=-\infty$.}
\ee
By (\ref{comp3}) and (\ref{comp4}), $\eta_\ve-\zeta$ has a positive maximum in $\Om_T$. Set $M=\sup_{\Om_T} (\eta_\ve-\zeta)$ and for $\al>0$, define
\ben
M_{\al}&=&\sup_{\Om\times \Om\times[0,T)} \left[ \eta_\ve(x,t)-\zeta(z,t)-\frac{ |x-z|^2}{2 \al}  \right].
\een
Then for each $\al>0$, $M_\al\ge M\ge \dl/2.$ Since both $\eta_\ve$ and $\zeta$ are bounded, there exists an $(x_\al,z_\al,t_\al)\in \overline{\Om}\times \overline{\Om}\times [0,T]$ which is a point of maximum, i.e.,
\eqRef{comp5}
M_\al=\eta_\ve(x_{\al},t_\al)-\zeta(z_\al,t_\al)-\frac{ |x_\al-z_\al|^2} {2 \al}
=\eta(x_{\al},t_\al)-\zeta(z_\al,t_\al)-\frac{ |x_\al-z_\al|^2} {2 \al}-\frac{\ve}{T-t_{\al}}.
\ee
It is known that as $\al\rightarrow 0$, $M_\al\rightarrow M$ and $|x_\al-z_\al|^2/\al\rightarrow 0.$ Let $(p, \s)\in \overline{\Om}_T$ be such that $x_\al,\;z_\al\rightarrow p$ and 
$t_\al\rightarrow \s$ (working with subsequences if needed) where $M=(\eta_\ve-\zeta)(p,\s).$ By (\ref{comp3}) and (\ref{comp4}), $(p,\s)$ is in $\Om\times[0,T]$ and hence, we may assume
that the sequences $x_\al,\;z_\al$ lie in a set compactly contained in $\Om$.
\vsp
\NI We now show that there are $0<t_0<T_0<T$ and $\al_0$, small enough, such that 
\eqRef{comp6}
t_0\le t_\al\le T_0,\;\;\mbox{for $0<\al\le \al_0$.}
\ee
Suppose that the first inequality is false. There is a sequence $\al_m\rightarrow 0$ such that $t_{\al_m}\rightarrow 0.$
Since $|x_{\al_m}-z_{\al_m}|^2/\al_m \rightarrow 0$, $x_{\al_m}$ and $z_{\al_m}\rightarrow q$, for some $q\in \overline{\Om}$ (work with subsequences if needed), we observe that 
\ben
0<\frac{\dl}{2}\le M &=& \lim_{\al_m \rightarrow 0} M_{\al_m}=\lim_{\al_m\rightarrow 0}  \left( \eta_\ve(x_{\al_m}, t_{\al_m})-\zeta(z_{\al_m},t_{\al_m})-\frac{|x_{\al_m}-z_{\al_m}|^2}{2\al} \right)\\
& \le & \limsup_{\al_m\rightarrow 0} \left ( \eta_\ve(x_{\al_m}, t_{\al_m})-\zeta(z_{\al_m}, t_{\al_m})-\frac{ |x_{\al_m}-z_{\al_m} |^2}{2 \al_m} \right) \\
&\le&(\eta_\ve-\zeta)(q,0)\le -\frac{\ve}{T}.
\een
The last inequality follows from (\ref{comp4}) and the hypothesis that $\eta\le \zeta$ on $P_T$. This is a contradiction and the assertion in (\ref{comp6}) holds. Next, since 
$\eta-\zeta$ is bounded in $\Om_T$ and 
$(T-t)^{-1}\rightarrow \infty$ as $t\uparrow T$, there is
a $T_0<T$ such that $t_{\al}\le T_0$ for all $\al>0.$
\vsp
\NI Recalling (\ref{comp5}), there are numbers $a_\al,\;b_\al\in \IR$ and $X_\al,\;Y_\al\in S(n)$ such that (see (\ref{def4}))
$$\left(a_\al, \frac{x_\al-z_\al}{\al}, X_\al \right)\in \overline{\mathcal{P}}^{+}_{\Om_T} \eta_\ve(x_\al,t_\al)\;\;\mbox{and}\;\;
\left(b_\al, \frac{x_\al-z_\al}{\al}, Y_\al\right)\in \overline{\mathcal{P}}^{-}_{\Om_T} \zeta(z_\al,t_\al).$$
Moreover,
$$a_\al=b_\al\;\;\;\;\mbox{and}\;\;\; \frac{-3}{\al} \left(\begin{array}{lr} I & 0\\  0& I \end{array}\right) \le \left(\begin{array}{lr} X_\al & 0\\  0& -Y_\al \end{array}\right)\le \frac{3}{\al}\left(\begin{array}{lr} I & -I\\  -I & I \end{array}\right).  $$
Thus $X_\al\le Y_\al$ and 
\eqRef{comp7}
\frac{\langle X_\al (x_\al-z_\al),(x_\al-z_\al)\rangle}{\al^2}\le\frac{\langle Y_\al (x_\al-z_\al),(x_\al-z_\al)\rangle}{\al^2}.
\ee
Using (\ref{comp1}), (\ref{comp2}), (\ref{comp6}) and with (\ref{def0}), we have
\ben
&& \langle X_\al (x_\al-z_\al),(x_\al-z_\al) \rangle \al^{-2} + \left |x_\al-z_\al \right| \al^{-4} - 3a_\al-3\ve (T-t_\al)^{-2} \ge 0, \\
\mbox{and}\;\;&&\langle Y_\al (x_\al-z_\al), (x_\al-z_\al) \rangle \al^{-2}+\left |x_\al-z_\al \right| \al^{-4}- 3a_\al \le 0.
\een
Using (\ref{comp7}),
\ben
3a_\al+\frac{3\ve}{(T-t_\al)^2}-\frac{|x_\al-z_\al|^4}{\al^4}\le 3a_\al-\frac{|x_\al-z_\al|^4}{\al^4}.
\een
This is a contradiction and $\eta\le \zeta$ in $\Om_T$. To complete the proof, we recall that $u=e^{\eta}$ and $v=e^{\zeta}$. The theorem holds. $\Box$
\vsp
\begin{rem}\label{compa}
The comparison principle in Theorem \ref{comp} continues to hold for equations of the type
$$\Df \eta+\s |D\eta|^4 \ge (\le)\; 3\eta_t,$$
where $-\infty<\s<\infty$ is a constant. $\Box$
\end{rem}
\begin{cor}\label{compc} Suppose that $\eta$ is upper semi-continuous and $\zeta$ is lower semi-continuous in $\Om_T\cup P_T$ and satisfy
$$\Df \eta+\s|D\eta|^4\ge 3 \eta_t\;\;\;\mbox{and}\;\;\;\Df \zeta+\s|D\zeta|^4\le 3\zeta_t\;\;\;\mbox{in $\Om_T$},$$
where $-\infty<\s<\infty$ is a constant.
Assume that $\eta$ and $\zeta$ are bounded in $\Om_T\cup P_T$ then
$$  \sup_{\Om_T}(\eta-\zeta)\le   \sup_{P_T}(\eta-\zeta).$$
Moreover, if $\eta$ and $\zeta$ are solutions then $\eta,\;\zeta\in C(\Om_T\cup P_T)$ and 
$$\sup_{P_T} |\eta-\zeta|=\sup_{\Om_T}|\eta-\zeta|.$$
In particular, if $\eta=\zeta$ on $P_T$ then $\eta=\zeta$ in $\Om_T$.
\end{cor}
\NI{\bf Proof:} Let $k=\sup_{P_T}(\eta-\zeta)$ and $\zeta_k=\zeta+k$. Then $\zeta_k$ continues to be a super-solution and $\zeta_k\ge \eta$ on $P_T$. By Theorem \ref{comp} and Remark \ref{compa}, $\sup_{\Om_T}(\eta-\zeta_k)\le \sup_{P_t}(\eta-\zeta_k)=0,$ and the claim holds. The second claim now follows quite easily. $\Box$

\begin{rem}\label{pbcmp}
Let $u$ and $v$ be as in Theorem \ref{comp}, i.e., $\Df u\ge 3u^2 u_t,\;\mbox{and}\;\Df v\le 3v^2 v_t,\;\mbox{in $\Om_t$},$
and $u,\;v$ be positive in $\Om_T\cup P_T$. Then
$$\sup_{P_T}\frac{u}{v}\ge\sup_{\Om_T}\frac{u}{v} \;\;\;\mbox{or}\;\;\;\sup_{P_T}\frac{u-v}{v}\ge \sup_{\Om_T}\frac{u-v}{v}.$$
This follows from Corollary \ref{compc} by writing $\eta=\log u$ and $\zeta=\log v$. For solutions, we obtain
$$\sup_{P_T} \frac{u}{v}=\sup_{\Om_T}\frac{u}{v}\;\;\;\mbox{and}\;\;\;\inf_{P_T}\frac{u}{v}=\inf_{\Om_T}\frac{u}{v}.$$
This leads to uniqueness for positive solutions to (\ref{pb1}). 
\vsp
\NI The above result is analogous to the result in Lemma 2.3 in \cite{BL}, where we showed a comparison principle for positive solutions of problems related to the eigenvalue problem for the infinity-Laplacian.  
$\Box$
\end{rem}

\begin{rem}\label{compz}
The comparison principle in Theorem \ref{comp} can be extended to include the case when $\inf_{P_T} v= 0.$ 
We assume that both $u$ and $v$ are positive in $\Om_T$. Let 
$$N=\{(y,s)\in P_T:\;v(y,s)=0\}.$$
Suppose that 
\eqRef{compz1}
\sup_{(y,s)\in P_T\setminus N}\frac{u(y,s)}{v(y,s)}\le 1,\;\mbox{and}\; \limsup_{(x,t)\rightarrow (y,s)}\frac{u(x,t)}{v(x,t)}\le 1,\;\;\mbox{$\forall(y,s)\in N$ and $(x,t)\in \Om_T$.}
\ee
Then $u\le v$ in $\Om_T$. 
\vsp
\NI To see this, choose $\bar{T}$ close to $T$ with $\bar{T}<T$; fix $\bar{T}.$ For every $\dl>0$, small, set $\Om^\dl=\{x\in \Om:\;\mbox{dist}(x,\p\Om)>\dl\}$ and $U_\dl=\Om^\dl\times(\dl,\bar{T}-\dl)$. Call $H_\dl$ the parabolic boundary of $U_\dl$. For a small and fixed $\ve>0$, set $v_\ve=(1+\ve)v$. Remark \ref{pbcmp} implies
$$\sup_{H_\dl}\frac{u}{v_\ve}\ge\sup_{U_\dl}\frac{u}{v_\ve}.$$ 
The right hand side increases as $\dl$ decreases to zero. Since $u>0$ is upper semi-continuous and $v>0$ is lower semi-continuous in $\Om_T\cup P_T$, $u/v$ is upper semi-continuous in $\Om_T\cup (P_T\setminus N)$. We claim that for all
$0<\dl\le \dl_0$, $\dl_0$ small enough,
$$\sup_{U_\dl}\frac{u}{v_\ve}\le \frac{1}{1+(\ve/2)}.$$
If not, there is a sequence $\dl_k\downarrow 0$ as $k\rightarrow \infty$ such that $(x_k,t_k)\in \overline{U}_{\dl_k}$ and $(u/v_\ve)(x_k,t_k)>(1+(\ve/2))^{-1}.$ 
For $k=1,2,\cdots$, select $0<\nu_k<\dl_k$, then $(x_k,t_k)\in U_{\nu_k}$. Invoking Remark \ref{pbcmp}, there is a point $(y_k,s_k)\in H_{\nu_k}$ 
such that $(u/v_\ve)(y_k,s_k)>(1+(\ve/2))^{-1}.$ Note that $\overline{U}_{\dl_k}\subset U_{\nu_k}\subset \overline{U}_{\nu_k} \subset\Om_{\bar{T}}.$ Also, we may select $\nu_k\downarrow 0.$
\vsp
\NI Let $(y,s)\in P_T$ be a limit point of $(y_k,s_k)$, then $s\le \bar{T}$. If $(y,s)\in N$, working with a subsequence of $(y_k,s_k)$, we obtain
$$\frac{1}{1+(\ve/2)}\le \limsup_{(y_k,s_k)\rightarrow (y,s)}\frac{u(x,t)}{v_\ve(x,t)}\le \frac{1}{1+\ve}.$$
If $(y,s)\in P_T\setminus N$ then $u/v$ is upper semi-continuous at $(y,s)$ and 
$$\frac{1}{1+(\ve/2)}\le \limsup_{(y_k,s_k)\rightarrow (y,s)}\frac{u(x,t)}{v_\ve(x,t)}\le \frac{u(y,s)}{v_\ve(y,s)}\le \frac{1}{1+\ve}.$$
Both cases lead to contradiction and our claim holds. Since $\ve$ is arbitrary,
\eqRef{compz1}
\lim_{\dl\downarrow 0} \left( \sup_{U_\dl} \frac{u}{v} \right)\le 1,
\ee
Since $\overline{U}_\dl\subset \Om_{\bar{T}}$, we have $u\le v$ in $\Om_{\bar{T}}.$
Suppose that, for some $(y,s)\in \Om_T$, we have $(u/v)(y,s)>1$. Clearly, $s<T$, selecting $s<\bar{T}<T$, \eqref{compz1} leads to a contradiction. The conclusion holds. $\Box$
\end{rem}
\vsp
\NI We end this section by proving a weaker version of the strong minimum principle, see, for instance, \cite{LN}. We show that if a super-solution attains its minimum, on $\Om_T$, at some point $y\in \Om$ and at some instant $s$, then, at $y$, it has the same value at all instants prior to $s$. This will be shown by using the comparison principle in Lemma \ref{comp}. 
Recall (\ref{def3}) for the definition of $F_{r,\tau}(x,t)=B_r(x)\times(t,t+\tau)$. This lemma uses the minimum on the parabolic boundary $P_T$. One could easily state the same with the minimum being taken on $\Om_T$.

\begin{lem}\label{minm}
Let $\zeta$ be lower semi-continuous in $\Om_T\cup P_T$ and satisfy
$$\Df \zeta+\s|D\zeta|^4\le 3 \zeta_t\;\;\;\mbox{in $\Om_T$},$$
where $-\infty<\s<\infty$ is a constant.
Let $m=\inf_{P_T} \zeta>-\infty$. If there is a point $(y,s)\in \Om_T$ such that
$\zeta(y,s)=m$ then 
$$\zeta(y,t)=m,\;\;0< t\le s.$$

\end{lem}

\NI{\bf Proof:} The function $\zeta_m=\zeta-m$ in $\Om_T$ satisfies
$$\Df \zeta_m +\s|D\zeta_m|^4\le 3(\zeta_m)_t, \;\;\mbox{in $\Om_T$}.$$
By Remark \ref{compa}, $\zeta_m\ge 0$ in $\Om_T$. We will show that given any $0<\tau<s$, $\zeta_m(y,t)=0,\;\forall t\in[\tau,s]$. Fix $\tau$, and  assume to the contrary. Then there exist $0<\ve\le s-\tau$ and $\dl>0$ such that $\zeta_m(y,s-\ve)> 2\dl$.
By lower semi-continuity, there is a $\rho>0$,small, such that 
\eqRef{minm1}
\zeta_m(x,s-\ve)\ge \dl,\;\;\forall x\in B_{\rho}(y).
\ee
Next, define
\eqRef{minm2}
\psi(x,t)=K  \left(  \rho^2-|x-y|^2  \right)^2 h(t),\;\;\forall(x,t)\in F_{\rho, 4\ve/3}(y,s-\ve),
\ee
where
\eqRef{minm3}
h(t)=\left(1-\frac{t-s+\ve}{2\ve}\right),\;\;\mbox{and}\;\;0<K\le \min \left( \sqrt{ \frac{3}{128\ve\rho^4(1+4|\s|\rho^4)} },\;\;\frac{\dl}{\rho^4} ,\;1 \right).
\ee
We observe the following:
\bea\label{minm4}
&&\psi(y,s)=\frac{K\rho^4}{2}\le \frac{\dl}{2},\;\;\sup_{x\in B_\rho(y)}\psi(x,s-\ve)\le K\rho^4\le \dl,   \nonumber\\
\mbox{and}&&\psi(x,t)=0,\;\forall(x,t)\in \p B_{\rho}(y)\times [s-\ve,\; s+\ve/3].
\eea
Set $r=|x-y|$ and write $\psi(x)=\psi(r)$ in (\ref{minm2}), and then use (\ref{minm3}) to calculate
\ben
\Df \psi + \s|D\psi|^4-3\psi_t  &\ge&  (\psi_r)^2\psi_{rr}-|\s| |\psi_r|^4-3\psi_t\\
&\ge&  \frac{3K}{2\ve} \left( \rho^2-r^2 \right)^2 -64K^3 r^2\left(\rho^2-r^2 \right)^2 (\rho^2-3 r^2) h(t)^3\\
& -& |\s| (4K)^4 r^4(\rho^2-r^2)^4 h(t)^4.
\een
Rearranging and noting that $1/3\le h(t)\le 1,\;\forall t\in[s-\ve, s+\ve/3]$,
\ben
\Df \psi + \s|D\psi|^4-3\psi_t 
\ge K(\rho^2-r^2)^2 \left( \frac{3}{2\ve}- 64 K^2 r^2 (\rho^2-3r^2) -4^4|\s| K^3 r^4(\rho^2-r^2)^2 \right).
\een
Using $r^2(\rho^2-3r^2)\le \rho^4,\;0\le r\le \rho$ and $K\le 1$, we find
\ben
\Df \psi + \s|D\psi|^4-3\psi_t &\ge& K(\rho^2-r^2)^2 \left( \frac{3}{2\ve}-64K^2 \rho^4 - 4^4|\s| K^3 \rho^8 \right)\\
&\ge& K(\rho^2-r^2)^2 \left( \frac{3}{2\ve}-64K^2 \rho^4 \left(1+  4|\s| \rho^4 \right) \right).
\een
By (\ref{minm1}), (\ref{minm3}) and (\ref{minm4}), we have
$$
\Df \psi+\s|D\psi|^4\ge 3\psi_t,\;\;\mbox{in $F_{r,4\ve/3}(y,s-\ve)$,}$$
$\psi(y,s)>0,$ and $\psi(x,t)=0,\;\forall(x,t)\in B_{\rho}(y)\times[s-\ve, s+4\ve/3].$ Since $\zeta_m\ge \psi$ on the parabolic boundary of $F_{\rho,4\ve/3}(y,s-\ve)$, from Remark \ref{compa} we have that
$\zeta_m\ge \psi$ in $F_{\rho,4\ve/3}(y,s-\ve)$ and $\zeta(y,s)\ge \psi(y,s)>0.$ This is a contradiction and $\zeta_m(y,t)=0,\;0<\tau<t<s.$ Since $\tau$ is arbitrary the claim holds. $\Box$
\vsp
\begin{rem}\label{maxm}
An analogous result holds for sub-solutions. If $\eta$ is upper semi-continuous in $\Om_T\cup P_T$ and solves
$$\Df \eta+\s|D\eta|^4\ge 3 \eta_t,\;\;\;\mbox{in $\Om_T$}$$
and if there is a point $(y,s)\in \Om_T$ such that $\eta(y,s)=\sup_{\Om_T} \eta$ then $\eta(y,t)=\sup_{\Om_T}\eta$ for $0<t\le s.$
Set $M=\sup_{\Om_T}\eta$ and call $\eta_M=M-\eta$, then
$$\Df\eta_M-\s|D\eta_M|^4\le 3(\eta_M)_t\;\;\;\mbox{in $\Om_T$}$$
and $\eta_M\ge 0$. The claim follows from Lemma \ref{minm}. $\Box$
\end{rem}
\vsp
\begin{rem}\label{maxm1} Using Lemma \ref{pde}, the conclusions of Lemma \ref{minm} and Remark \ref{maxm} apply to sub-solutions and super-solutions of 
$\Df \phi= 3\phi^2 \phi_t$. $\Box$
\end{rem}
\vspp

\section{Asymptotic behaviour and some results on $\IR^n\times(0,T)$}
\vsp
\NI In this section, we prove some results when $\Om=\IR^n$, and also some asymptotic results on a bounded domain $\Om$. We do not address the question of existence of solutions in the two situations we discuss below.
\vsp
\NI{\bf (a) Unbounded domain $\Om_T=\IR^n\times (0,T)$}
\vsp
\NI We take $\Om=\IR^n$. We continue to use the notation $\Om_T=\IR^n\times(0,T).$ In the event $T=\infty$, we write
$\Om_\infty=\Om\times(0,\infty).$
\vsp
\NI To prove our results, we make use of Theorem \ref{rad1} and Lemmas \ref{bnd}, \ref{rad2}, \ref{rad4} and \ref{sol1}. We first present some simple examples.
Suppose that $\psi_\lam(x)>0$, radial, solves the problem
$$\Df \psi_\lam-\lam \psi_\lam^3=0,\;\;\mbox{in $\Om$, and $\psi(0)=\dl$.}$$ 
where $\lam\ge 0$ and $\dl\ge 0$. Also, the solution $\psi_\lam$ can be extended to all of $\IR^n$.
Then the function $\phi(x,t)=\psi_\lam(x) e^{-\lam t/3}$ solves
$$\Df \phi=3\phi^2 \phi_t,\;\;\;\mbox{in $\Om_T$}.$$
See Remark \ref{radp} and Lemmas \ref{rad4} and \ref{sol1}. Also shown in Remark \ref{radp} is that
this solution is unbounded in $\Om_T$. Suppose $M>0$ is a constant,  then $\forall k>0$, $\phi_+(x,t)=M e^{kt}$ is a super-solution, and $\phi_-=Me^{-kt}$ is a sub-solution. Both are bounded in $\Om_T$.
\vsp
\NI We begin by showing that if the initial data is bounded then all bounded non-negative solutions satisfy a maximum principle. A slight generalization will be proven later. 
Let $y\in \Om$; for $r>0$ and $T>0$, call $B_{r,T}(y)=B_r(y)\times (0,T)$ and let $P_{r,T}(y)$ denote the parabolic boundary of $B_{r,T}(y)$. We use Lemma \ref{rad4} and Remark \ref{radp} for our proof.

\begin{lem}\label{unbd0}
Let $\Om=\IR^n$ and $T>0.$ Suppose that $\phi\in C(\Om_T\cup (\Om\times\{0\}) ),\;\phi> 0,$ solves 
$$\Df \phi=3\phi^2 \phi_t,\;\;\mbox{in $\Om_T$ and $\phi(x,0)=f(x),\;\forall x\in \Om$,}$$
where $f:\Om\rightarrow \IR^+$ is continuous and bounded. If $\phi$ is bounded then
$$\inf_{\Om} f\le \phi(x,t)\le \sup_\Om f, \;\forall(x,t)\in \Om_T.$$
\end{lem}
\NI{\bf Proof:} Set $\nu=\inf_{\Om} f,\;\mu=\sup_{\Om} f$ and $\sup_{\Om_T} \phi=M$. For a fixed $y\in \Om$, set 
$r=|x-y|$ and, for $\lam>0$, define $\psi_\lam(x)=\psi_\lam(r),\;\forall x\in \Om,$ by
\eqRef{unbd1}
\int_{\mu}^{\psi_\lam(r)} \frac{ds}{ (s^4-\mu^4)^{1/4} }=\lam^{1/4} r,\;\;\forall r\ge 0.
\ee
Then
$$\Df \psi_\lam-\lam \psi_\lam^3=0,\;\; \mbox{ in $\Om,$ and $\psi_\lam(0)=\mu$.}$$
Moreover, $\psi_\lam(r)$ is continuous and increasing in $r$. These follow from Lemma \ref{rad4}. From Remark \ref{radp}, for each $\lam>0$, there is an $R=R(\lam)>0$ such that
$\psi_\lam(R)=\max(M+1,\;3\mu).$ Define
$$\Gamma_\lam(x,t)=\psi_\lam(x) e^{\lam t/3}=\psi(r)e^{\lam t/3},\;\;\forall(x,t)\in \Om_T.$$
By Lemma \ref{sol1}, we have 
$$
\Df \Gamma_\lam=3\Gamma_\lam^2(\Gamma_\lam)_t,\;\mbox{in $B_{R,T}(y)$},
$$
and $\Gamma_\lam\ge \phi$ on $P_{R,T}(y)$. By Theorem \ref{comp}, $\Gamma_\lam\ge \phi$, in $B_{R,T}(y),$ and in particular, 
$\phi(y,t)\le \psi_\lam(0)e^{\lam t/3}=\mu e^{\lam t/3}.$ By Remark \ref{radp} or by (\ref{unbd1}), for every $\lam>0$, 
if $R$ is such that $\psi_\lam(R)=\max(M+1,3 \mu)$ then $\lam\rightarrow 0$ as $R\rightarrow \infty$. Thus, we obtain
$\phi(y,t)\le \mu,\;\forall (y,t)\in \Om_T.$ 
\vsp
\NI To show the lower bound, we use Theorems \ref{eignv} and \ref{rad1}. Assume that $\nu>0$, and let $y\in \Om$, $r=|x-y|$ and $R>0$, large. Let $\lam_R>0$ be the first eigenvalue of $\Df$ on $B_R(y)$ and $\xi_R(r)>0,\;$ a first radial eigenfunction, with $\xi_R(R)=0$, i.e., $\Df \xi_R+\lam_R \xi^3=0$, in $B_R(y)$. 
Scale $\xi_R(0)=\nu.$ By Lemma \ref{sol1}, the function
\eqRef{unbd11}
\Theta_R(x,t)=\xi_R(x)e^{-\lam_R t/3}=\xi(r) e^{-\lam_R t/3},
\ee
solves
\eqRef{unbd12}
\Df \Theta_R=3\Theta_R^2 (\Theta_R)_t,\;\mbox{in $B_{R,T}(y),$ $\Theta_R(R,t)=0,\;\forall 0\le t\le T,$}
\ee
 and $\Theta_R(x,0)\le\nu,\;\forall x\in B_R(y).$ Thus, $\Theta_R\le \phi$, on $P_{R,T}(y),$ and using Remark \ref{compz}, we have $\Theta_R\le \phi$ in $B_{R,T}$. Taking $r=0$, 
we get $\phi(y,t)\ge \nu e^{-\lam t/3},\;\forall\;0<t<T$. The relation between $\lam_R$ and $R$ (see Theorem \ref{rad1})
$\displaystyle{ \lam_R^{1/4}R=\int^1_0 (1-s^4)^{-1/4}\;ds,} $
implies that $\lam_R\rightarrow 0$ as $R\rightarrow \infty$. Thus, the lower bound in the lemma holds. $\Box$   
\vsp
\NI In the next lemma, we replace the restriction of boundedness by an exponential type growth and obtain the maximum principle. 

\begin{lem}\label{unbd3}
Let $T>0$ and $\Om_T=\IR^n\times(0,T).$ Suppose that $\phi\in C(\Om_T\cup (\Om\times\{0\}) ),\;\phi> 0,$ solves 
$$\Df \phi=3\phi^2 \phi_t,\;\;\mbox{in $\Om_T$ with $\phi(x,0)=f(x),\;\forall x\in \Om$,}$$
where $f:\Om\rightarrow \IR^+$ is continuous and bounded. Suppose that $h:\IR^+\times(0,T)\rightarrow \IR^+$ is such that 
$$\frac{\sup_{0<t<T}h(r,t)}{r}=o(1)\;\;\mbox{as $r\rightarrow \infty$.}$$
If $\phi(x,t)\le e^{h(|x|,t)},\;\;\forall(x,t)\in \Om_T,$ then
$$\inf_{\Om} f\le \phi(x,t)\le \sup_\Om f,\;\forall (x,t) \in \Om_T.$$
\end{lem}
\NI{\bf Proof:} We assume that $h$ is not bounded. Set $\nu=\inf_{\Om} f$, $\mu=\sup_{\Om} f$ and assume that $\nu>0.$ Let $y\in \Om$. For any $r>0$ and $T>0$, 
call $B_{r,T}(y)=B_r(y)\times(0,T)$ and $P_{r,T}(y)$ its parabolic boundary.
\vsp
\NI In order to show that $\phi\ge \nu$, we use (\ref{unbd11}) and (\ref{unbd12}), see Lemma \ref{unbd0}. Let $y\in \Om$ and $R>0$. Since $\phi\ge 0$ on $\p B_R(y)\times[0,T)$, the function $\Theta_R$, as defined in (\ref{unbd11}), lies below $\phi$ in $B_R(y)\times[0,T).$ Rest of the proof for the lower bound is the same as in Lemma \ref{unbd0}. 
\vsp
\NI For the upper bound, we use Remark \ref{radp}. Set 
$\al(|x|)=\sup_{0<t<T} h(|x|,t)$. Let $\ve>0$ be small. Choose $R$ large, so that  $\al(|x|)\le \ve |x|,\;\;\forall |x|\ge R$. This implies that
\eqRef{unbd4}
\nu\le \phi(x,t)\le e^{\al(|x|)}\le e^{\ve |x|},\;\;\forall |x|\ge R.
\ee
Let $y\in \Om$; set $r=|x-y|,\;x\in \Om. $ Call $\rho=R-|y|$.
By Lemma \ref{rad4} and Remark \ref{radp}, there is a $\lam>0$ and a function $\psi_\lam(x)=\psi_\lam(r)>0$ such that 
\eqRef{unbd5}
\Df \psi_\lam-\lam \psi_\lam^3=0,\;\;\mbox{in $B_\rho(y)$, $\psi_\lam(0)=\mu$ and $\psi_\lam(\rho)= e^{\ve R} $}.
\ee
Note that $B_\rho(y)\subset B_R(o)$. We extend $\psi_\lam$ to all of $[0,\infty)$ as a solution of the equation in (\ref{unbd5}). Note that $\psi_\lam(r)$ is increasing in $r$.
By Lemma \ref{sol1}, the function
$$\Gamma_\lam(x,t)=\psi_\lam(x) e^{\lam t/3}$$
solves
$$\Df \Gamma_\lam=3\Gamma_\lam^2 (\Gamma_\lam)_t,\;\mbox{in $B_{\rho,T}(y),$ and $\Gamma_\lam\ge \phi$ on $P_{\rho,T}(y)$.}$$
Applying Theorem \ref{comp}, $\Gamma_\lam\ge \phi$ in $B_{\rho,T}(y)$.
Using Remark \ref{radp}, (\ref{rad13}) and (\ref{unbd5}) there is a constant $K>0$ independent of $\rho$ (for $R$ large enough), such that 
$$K\mu e^{\lam^{1/4} \rho/\s}\le \psi_\lam(\rho)\le e^{\ve \rho+|y|},\;\;\;\mbox{where}\;\;\;\s=\frac{2}{(2^4-1)^{1/4}}.$$
Thus,
$$\lam^{1/4}\le \s \left( \ve\left(1+\frac{|y|}{\rho}\right)-\frac{\log(K\mu)}{\rho}    \right) \le 2\s \ve,$$
if $\rho$ is chosen large enough. We conclude that 
$$\phi(y,t)\le \Gamma_\lam(y,t)=\Gamma_\lam(0,t)=\psi_\lam(y)e^{\lam t/3}\le \mu e^{(2\s \ve)^4 T/3}.$$ 
Since $\ve$ is arbitrary, $\phi(y,t)\le \mu.\;\;\;\Box$
\vsp
\NI{\bf (b) Asymptotic behavior $\Om_\infty=\Om\times(0,\infty).$}
\vsp
\NI Our next task is to address the issue of asymptotic behavior of solutions. Let $\Om\subset \IR^n$ be a bounded domain. 
We discuss $\lim_{t\rightarrow \infty} \phi(x,t)$ where $\phi$ solves the problem in (\ref{pb1}) and where $f$ and $g$ satisfy certain conditions. Call $\Om_{\infty}=\Om\times[0,\infty)$ and 
$P_\infty=(\p\Om\times [0,\infty))\cap (\Om\times \{0\}),$ and set 
\eqRef{asym0}
m=\min( \inf_{\overline{\Om}} f, \inf_{\p\Om\times (0,\infty)} g)\;\;\;\;\mbox{and}\;\;\;\;\;M=\max(\sup_{\overline{\Om}} f, \sup_{\p\Om\times (0,\infty)} g ).
\ee
Our first result will provide information about $\phi$ for large $t$, when $g=0$ on $\p\Om\times[0,\infty)$, see Lemma \ref{asym}. Later we prove a more general result.
We begin with an estimate that will be used in Lemma \ref{asym}.
\begin{lem}\label{asym01}
Let $\Om\subset \IR^n$ be a bounded domain and $\psi\in C(\overline{\Om}),\;\psi>0$ solve
$$\Df \psi+\lam \psi^3=0,\;\;\mbox{in $\Om,$ and $\psi=\dl$ on $\p\Om$,}$$
where $\lam>0$ and $\dl\ge 0.$ Then there is a constant $C>0$, depending only on $\Om$, such that
for any $x\in \Om$, close to $\p\Om$, 
$$\psi(x)-\dl \le C(\sup_\Om \psi)\mbox{dist}(x,\p\Om).$$
\end{lem}
\NI {\bf Proof:} We will utilize the comparison principle stated in Lemma \ref{comp1}, also see \cite{BMO2}. Suppose that $L=\sup_\Om \psi$. Let $z\in \p\Om$; 
set $r=|x-z|,\;\forall x\in \IR^n,$ and
$R_z=\sup_{x\in \Om} |x-z|$. Then $\Om\subset B_{R_z}(z)$.
Call $\s=3^4/4^3$, define $u(x)=u(r)$ in $0\le r\le R$, and set $\forall x\in B_R(z)$,
\eqRef{asym011}
u_z(r)=\dl+K_z \left( R_z^{4/3}-(R_z-r)^{4/3} \right),\;\;\mbox{where}\;K_z= L \left( \max\left( R_z^{-4/3},\; (\lam \s)^{1/3}  \right)  \right).
\ee
Observe that there are $0<R_0<R_1<\infty$ such that $R_0\le R_z\le R_1,\;\forall z\in \p\Om$.
Thus, $0<K_z\le K_0$, for some $K_0<\infty,$ where $K_0=L \left( \max\left( R_0^{-4/3},\; (\lam \s)^{1/3}  \right)  \right).$ Note that $L$ depends on $\lam$, $\dl$ and $\Om$, see Lemma \ref{bnd}. 
\vsp
\NI First observe that $u_z(x)\ge \dl,\;\forall x\in B_{R_z}(z)$. Next, differentiating $u_z$, as in (\ref{asym011}), leads to
\eqRef{asym02}
\Df u_z=\left(\frac{du_z}{dr}\right)^2 \frac{d^2u_z}{dr^2}\le -\lam L^3,\;\mbox{in $B_{R_z}(z)\setminus \{z\}$, and $u_z(0)=u_z(z)=\dl$, and $u_z(R_z)\ge L$.}
\ee
We recall that $\Df \psi=-\lam \psi^3\ge -\lam L^3$, in $\Om$, and $\psi=\dl$ on $\p\Om$. Applying Lemma \ref{comp1} in 
$\Om$, $\psi\le u_z$, in $\Om$, and recalling that for a fixed $z\in\p\Om$ that $r=|x-z|$, we get
$$\psi(x)-\dl\le K_0\left( R_z^{4/3}- (R_z- r)^{4/3} \right)\le \left(\frac{4K_0R_1^{1/3}}{3}\right) |x-z|,\;\;\forall x\in \Om.$$
Thus, for any $x\in \Om$, close to $\p\Om$, $\psi(x)-\dl\le CL \mbox{dist}(x,\p\Om)$, where 
\eqRef{asym03}
C=\frac{4R_1^{1/3}}{3}\left( \max\left( R_0^{-4/3},\; (\lamo \s)^{1/3}  \right)  \right). \;\;\;\Box
\ee
\vsp
\begin{lem}\label{asym}
Let $\Om\subset \IR^n$ be a bounded domain and $f:\overline{\Om}\rightarrow \IR^+$ be continuous. Suppose that $\phi\in C(\Om_\infty\cup P_\infty),\;\phi>0,$ solves the following problem:
$$\Df \phi\ge3\phi^2 \phi_t,\;\mbox{in $\Om_\infty$,\;$\phi(x,0)=f(x),\;\forall x\in \overline{\Om},$ and $\phi(x,t)=0,\;\;\forall(x,t)\in \p\Om\times(0,\infty)$.}$$
Assume that $m$ and $M$, as in (\ref{asym0}), satisfy $0<m \le M < \infty$. Then $\sup_\Om\phi(x,t)$ is decreasing in $t$,
$$\lim_{t\rightarrow \infty} (\sup_{\Om}\phi(x,t))=0,\;\;\;\mbox{and}\;\;\;
\lim_{t\rightarrow \infty} \frac{\log(\sup_\Om \phi(x,t) )}{t} \le -\frac{\lam_\Om}{3}. $$
\end{lem}

\NI{\bf Proof:} For every $T>0$, working in $\Om_T$, the maximum principle in Lemma \ref{max} implies that $m\le \phi(x,t)\le M,\;\forall (x,t)\in \Om_T.$ For $t\ge 0$, define $\mu_t=\sup_{x\in \overline{\Om}}\phi(x,t)$, and for any $0\le T<\bar{T}<\infty$, consider the domain
$\Om\times (T, \bar{T})$. Since $\phi(x,t)=0,\;\forall(x,t)\in \p\Om\times[0,\infty)$, Lemma \ref{max} implies that for $T\le \tau\le \bar{T}$,
$$\mu_\tau\le \sup_{\Om\times(T,\bar{T})}\phi(x,t)=\sup_{\p\Om\times[T,\bar{T})}\phi(x,t)=\sup_{\overline{\Om}}\phi(x,T)=\mu_T.$$ 
Hence, $\mu_t$ is decreasing in $t$.
\vsp
\NI We now show the remaining assertions. For any $\rho>0$, small, call $\Om_\rho=\{x\in \Om:\;\mbox{dist}(x,\p\Om)> \rho \}$. Fix $T>0$ and select $0<\ve\le \mu_T/4$, small.
Since $\phi(x,T)$ is continuous on $\overline{\Om}$ and $\phi(x,T)=0,\;\forall x\in \p\Om$, there is a $\rho>0$, small, such that $\phi(x,T)\le \ve,\;\forall x \in\overline{\Om}\setminus \Om_{\rho}$. This follows from compactness.
By Lemma \ref{exst}, for every $0<\lam<\lam_\Om$, there is a unique $\psi_\lam \in C(\overline{\Om}),\;\psi_\lam>0,$ that solves
\eqRef{asym1}
\Df \psi_\lam+\lam \psi_\lam^3=0,\;\;\mbox{in $\Om$, with $\psi_\lam=\ve$ on $\p\Om$.}
\ee
Since $\psi_\lam>\ve$ in $\Om$, $\phi(x,T)\le \psi_\lam(x),\;\;x\in \Om\setminus \Om_\rho.$ 
\vsp
\NI Our goal is to choose $\lam<\lam_\Om$ such that $\psi_\lam\ge \phi(x,T),\;\forall x\in \overline{\Om}$. To see this, we use the estimate in Lemma \ref{bnd} and select 
$\lam$ so that $\sup_{\Om} \psi_\lam>\mu_T$. Next, we apply the estimate in Lemma \ref{asym01} to $\psi_\lam$. If $\mbox{dist}(x,\p\Om)\le 1/(3C)$ 
(where $C$ is as in (\ref{asym03})) then 
$$\psi_\lam(x)\le \ve+\frac{\sup_{\Om} \psi_\lam}{3}\le \frac{\mu_T}{4}+\frac{\sup_{\Om} \psi_\lam}{3}<\sup_\Om \psi_\lam,\;\;\forall x\in \overline{\Om}\setminus \Om_{1/(3C)}.$$ 
This means, $\psi_\lam$ attains its maximum
at least $1/(3C)$ away from the boundary. If we choose $\rho$ small enough, then $\psi_\lam$ attains its maximum in $\Om_\rho$. Fix $\rho$.
Since $\Df \psi_\lam\le 0$ in $\Om$, $\psi_\lam$ satisfies the Harnack inequality (see \cite{ACJ}) and
$$
\inf_{\Om_\rho} \psi_\lam\ge K \sup_{\Om}\psi_\lam\ge K \mu_T,
$$
for some $K>0$ depending only on $\rho.$ In (\ref{asym1}), select $\lam$ closer to $\lam_\Om$, if needed, to ensure that $\inf_{\Om_\rho} \psi_\lam\ge \mu_T.$ With this choice
of $\lam$, $\phi(x,T)\le \psi_\lam(x),\;\forall x\in \overline{\Om}.$ We fix $\lam$ for what is to follow.
Now consider the function
$$\Gamma(x,t)=\psi_\lam(x) e^{-\lam (t-T)/3},\;\;\mbox{in $\overline{\Om}\times [T,\infty)$.}$$
Note that $\Gamma(x,T)\ge \phi(x,T),\;\forall x\in \overline{\Om}$, $\Gamma(x,t)>0,\;\forall(x,t)\in \p\Om\times[T,\infty),$ and by Lemma \ref{sol1}, 
$\Df \Gamma=3\Gamma^2 \Gamma_t$, in $\Om\times(T,\infty)$. Since $\phi(x,t)=0,\;\forall (x,t)\in \p\Om\times[0,\infty)$, 
Theorem \ref{comp} and Remark \ref{compz} imply that $\Gamma(x,t)\ge \phi(x,t),\;\forall(x,t)\in \Om\times(T,\infty).$ Letting $t\rightarrow \infty$, 
we obtain $\lim_{t\rightarrow \infty} (\sup_\Om \phi(x,t))=0.$ Moreover,
$$\lim_{t\rightarrow \infty} \frac{\log(\sup_\Om \phi(x,t)  )}{t}\le  \lim_{t\rightarrow \infty} \frac{\log(\sup_\Om \Gamma(x,t)  )}{t}=  -\frac{\lam}{3}.$$
Working the above arguments with $\lam$ closer to $\lam_\Om$, we get
$$\lim_{t\rightarrow \infty} \frac{\log(\sup_\Om \phi(x,t)  ) }{t}\le -\frac{\lamo}{3}.$$
$\Box$

\begin{rem}\label{asyms}
The requirement that $\phi$ be a sub-solution is necessary in Lemma \ref{asym}. To see this let $\psi>0$ be a first eigenfunction of $\Df$ in $\Om$, i. e.,
$$\Df \psi+\lam_\Om \psi^3=0,\;\;\mbox{in $\Om$ and $\psi=0$ on $\p\Om$.}$$
Define $\phi(x,t)=\psi(x)$ then $\Df \phi=-\lam \phi^3\le 3\phi^2 \phi_t=0.$ Clearly, $\sup_\Om \phi(x,t)$ does not decay to zero. $\Box$
\end{rem}

\NI We now state a generalization of Lemma \ref{asym}.

\begin{lem}\label{asymg}
Let $\Om\subset \IR^n$ be a bounded domain, $f:\overline{\Om}\rightarrow \IR^+$ and $g:\p\Om\times[0,\infty)\rightarrow \IR^+$ are continuous. Suppose that $\phi(x,t)\in C(\Om_\infty\cup P_\infty),\;\phi>0,$ solves
$$\Df \phi\ge 3\phi^2 \phi_t,\;\mbox{in $\Om_\infty$, $\phi(x,0)=f(x),\;\forall x\in \overline{\Om}$, and $\phi(x,t)=g(x,t),\;\forall(x,t)\in \p\Om\times [0,\infty)$}.$$
If $\lim_{t\rightarrow \infty} (\sup_{\p\Om} g(x,t))=0,$ then
$\lim_{t\rightarrow \infty} (\sup_{\Om\times [t,\infty)} \phi(x,t))=0.$
\end{lem}
\NI{\bf Proof:} By the maximum principle in Lemma \ref{max}, it is clear that $0<\phi\le M<\infty$, where $M=\max(\sup_{\Om} f,\;\sup_{\p\Om\times [0,\infty)}g)$, see (\ref{asym0}). By the hypothesis of the lemma, if $\ve>0$, small, then there is a $T_1>0$ such that 
\eqRef{asymg1}
0\le \sup_{x\in \p\Om} g(x,t)\le \frac{\ve}{2},\;\;\;\forall t\ge T_1.
\ee
We use the arguments in the proof of Lemma \ref{asym}. Since $\sup_{x\in\p\Om} \phi(x,T_1)<\ve$ (see (\ref{asymg1}))one can select $0<\lam<\lam_\Om$ (where $\lam=\lam(T_1)$) such that the solution $\psi_\lam>0$ of
\eqRef{asymg2}
\Df \psi_\lam+\lam \psi_\lam ^3=0,\;\;\mbox{in $\Om$, and $\psi_\lam(x)=\ve,\;\forall x\in \p\Om$,}
\ee
satisfies $\phi(x,T_1)\le \psi_\lam(x),\;\forall x\in \overline{\Om}$. Also, see Lemma \ref{asym01}. Also note that $\psi_\lam\in C(\overline{\Om})$ and $\psi_\lam(x)>\ve,\;\forall x\in \Om$. For $t\ge 0$, set $\mu_t=\sup_\Om \phi(x,t).$ 
\vsp
\NI Let $\bar{\lam}$ denote a value of $\lam$ such that the above holds (see (\ref{asymg2})) and $\psi$ denote the corresponding solution $\psi_{\bar{\lam}}$. We record that $0<\lamb<\lamo$ and $\phi(x,T_1)\le \psi(x,T_1),\;\forall x\in \overline{\Om}$. Define
\eqRef{asymg3}
g_1(t)= \frac{1}{2}\left(1+ \frac{  e^{\lamb(T_2-t)/3}-1 }{ e^{\lamb(T_2-T_1)/3}-1 }\right) ,\;\;\;\forall t\in [T_1,\;T_2].
\ee
We choose $T_2$ large enough so that $T_2>1+T_1$ and $e^{\lamb(T_2-T_1)/3}\ge 2.$ This choice is so made that $\sup_{\Om} g(x,t) \le \ve/4,\; \forall t\ge T_2.$ We note that
\eqRef{asymg4}
\psi|_{\p\Om}=\ve,\;\;\;\frac{1}{2}\le g_1(t)\le 1,\;\forall t\in [T_1,\;T_2],\;\;g_1(T_1)=1\;\;\mbox{and}\;\;g_1(T_2)=\frac{1}{2}.
\ee
Call $I_1=\Om\times (T_1,T_2)$ and $J_1=\p\Om\times[T_1,T_2).$ Define
$$\phi_1(x,t)=\psi(x) g_1(t),\;\;\forall (x,t)\in \overline{I_1}.$$
Arguing as in Lemma \ref{sol1}, using (\ref{asymg3}) and the comment about $T_2$, for $\forall (x,t)\in I_1$,
\ben
\Df \phi_1-3\phi_1^2 (\phi_1)_t  & =&  -\lamb \psi^3 g_1(t)^2    \left( g_1(t) - \frac{e^{\lamb(T_2-t)/3} }{2(e^{\lamb(T_2-T_1)/3}-1 ) }  \right) \\
&=& -\frac{\lamb \psi^3 g_1(t)^2}{2}  \left( \frac{e^{\lamb(T_2-T_1)/3}-2 }{ e^{\lamb(T_2-T_1)/3}-1  } \right)\le 0
\een
By our construction, see (\ref{asymg2}), (\ref{asymg3}), (\ref{asymg4}), the definition of $\psi$ and (\ref{asymg1}), we note
$$\phi(x,T_1)\le \psi(x)=\phi_1(x,T_1),\;\forall x\in \overline{\Om},\;\mbox{and}\;\phi(x,t)\le \frac{\ve}{2}\le \phi_1(x,t)=\psi(x)g_1(t)\le \ve,\;\forall(x,t)\in  J_1.$$
Since $\phi_1$ is a super-solution in $I_1$, Theorem \ref{comp} and Remark \ref{compz}  yield $\phi(x,t)\le \phi_1(x,t),\;\forall(x,t)\in I_1$, i.e.,
\eqRef{asymg5}
\phi(x,t)\le \psi(x),\;\forall(x,t)\in \overline{I}_1,\;\;\mbox{and}\;\;\phi(x,T_2)\le \phi_1(x,T_2)=\frac{\psi(x)}{2},\;\;\forall x\in \Om.
\ee
In particular, $\mu_{t}\le \sup_{\Om}\psi,\;\forall t\in [T_1,\;T_2],$ and $\mu_{T_2}\le \sup_{\Om} \psi/2.$ 
\vsp
\NI We now use induction. For $k=2,3,\cdots,$ call (see (\ref{asymg5}))
$$I_k=\Om\times(T_k, T_{k+1})\;\;\mbox{and}\;\;J_k=\p\Om\times[T_k, T_{k+1}).$$
For any $k$, once a choice of $T_k$ has been made with 
\eqRef{asymg51}
\sup_{x\in \p\Om}g(x,t)\le \frac{\ve}{2^k},\;\forall t\ge T_k,
\ee
the choice of $T_{k+1}$ is so made that $T_{k+1}>1+T_k$,
\eqRef{asymg6}
e^{\lamb(T_{k+1}-T_k)/3}\ge 2\;\;\;\mbox{and}\;\;\;\sup_{\p\Om} g(x,t)\le \frac{\ve}{2^{k+1}},\;\;\forall t\ge T_{k+1}.
\ee
Suppose that for some $k=1,2\cdots$, we have
\eqRef{asymg70}
\phi(x,T_k)\le \frac{\psi(x)}{2^{k-1}},\;\;\forall x\in \Om.
\ee
We show that 
\eqRef{asymg7}
\phi(x,t)\le \frac{\psi(x)}{2^{k-1}},\;\forall(x,t)\in \overline{I}_k, \;\;\;\mbox{and}\;\;\;\phi(x,T_{k+1})\le \frac{\psi(x) } {2^k},\;\forall x\in \Om.
\ee
Clearly, this will imply the conclusion of the lemma. To this end, define
$$g_k(t)=\frac{1}{2}\left(1+ \frac{e^{\lamb(T_{k+1}-t)/3}-1 }{e^{\lamb(T_{k+1}-T_k)/3}-1 } \right),\;\forall t\in[T_k,\;T_{k+1}].$$ 
Then
\eqRef{asymg8}
\frac{1}{2}\le g_k(t)\le 1,\;\;\forall t\in [T_k,\;T_{k+1}],\;\;g_k(T_k)=1\;\;\mbox{and}\;\;g_k(T_{k+1})=\frac{1}{2}.
\ee
Define
$$\phi_k(x,t)=\frac{\psi(x)g_k(t)}{2^{k-1}},\;\;\forall (x,t)\in \overline{I_k}.$$ 
Arguing as in Lemma \ref{sol1}, using the definition of $g_k(t)$ and (\ref{asymg6}), for $\forall(x,t)\in I_k$,
\ben
\Df \phi_k-3\phi_k^2 (\phi_k)_t  & = & -\frac{\lamb \psi^3 g_k(t)^2}{ 2^{3(k-1)} } \left( g_k(t) - \frac{e^{\lamb(T_{k+1}-t)/3} }{ 2(e^{\lamb(T_{k+1}-T_k)/3}-1) } \right)\\
&= & -\frac{\lamb \psi^3 g_k(t)^2}{ 2^{3(k-1)+1} }   \left( \frac{e^{\lamb(T_{k+1}-T_k)/3}-2 }{ e^{\lamb(T_{k+1}-T_k)/3}-1 } \right)\le 0.
\een
Thus, $\phi_k$ is super-solution in $I_k$. Using (\ref{asymg4}), (\ref{asymg51}) and (\ref{asymg8}),
$$\phi_k(x,T_k)=\frac{\psi(x)}{2^{k-1}}\ge \frac{\ve}{2^{k-1}},\;\forall x\in \overline{\Om},\;\;\mbox{and}
\;\;\frac{\ve}{2^k}\le \phi_k(x,t)=\frac{\psi(x) g_k(t)} {2^{k}}\le \frac{\ve}{2^{k-1}},\;\forall(x,t)\in J_k. $$
Using the inductive hypothesis (\ref{asymg70}) and recalling (\ref{asymg51}), $\phi_k\ge \phi$ on the parabolic boundary of $I_k$. By Theorem \ref{comp}, Remark \ref{compz} and (\ref{asymg8}),
$$\phi(x,t)\le \phi_k(x,t)=\frac{\psi(x) g_k(t)}{2^{k-1}}\le \frac{\psi(x)}{2^{k-1}},\;\forall(x,t)\in I_k,\;\;\mbox{and}\;\;\phi(x,T_{k+1})\le \frac{\psi(x)}{2^k},\;\forall x\in \Om.$$
Thus, (\ref{asymg7}) holds and by induction, for $k=2,3,\cdots,$ 
$$\phi(x,t)\le \frac{\psi(x)}{2^{k-1}},\;\;\forall (x,t)\in \Om \times [T_k,\;T_{k+1}].$$
The lemma holds. $\Box$.
\vsp

\section{Existence of positive solutions}
\vsp
\NI We now address the question of existence of positive solutions to (\ref{pb1}). Restated, let $\Om\subset \IR^n$ be bounded and $T>0,$ we seek a solution 
$\phi\in C(P_T\cup \Om_T),\;\phi>0,$ that solves
\eqRef{Exist0}
\Df \phi=3\phi^2\phi_t,\;\;\mbox{in $\Om_T$ and $\phi(x,t)=h(x,t),\;\forall(x,t)\in P_T$},
\ee
where we define $h(x,t)$ as follows. On $P_T$, set
$$h(x,t)= \left\{ \begin{array}{lcr} f(x), & \forall x\in \overline{\Om} \times\{0\},\\  g(x,t), & \forall(x,t)\in \p\Om\times [0,T). \end{array} \right. $$
By our notation, $h(x,0)=f(x),\;x\in \overline{\Om}.$ We assume that $h\in C(P_T)$. Set $m=\inf_{P_T}h$ and $M=\sup_{P_T}h$. In this section, we assume $0<m\le M<\infty$. 
Moreover, we will take $m<M$, for otherwise, $\phi(x,t)=m,\;\forall(x,t)\in \Om_T$ is the unique positive solution.
These notations introduced above will be followed in the rest of this section, in Appendix I and Appendix II, see Sections 6 and 7.
\vsp
\NI We recall the definitions of $\aleph_{sub}$ and $\aleph_{sup}$, see (\ref{Exist3}), (\ref{sub12}) and (\ref{sup11}) in Appendix I, see Section 6. 
We introduce the following notation. For every $(y,s)\in P_T$ and every $\ve>0$, the function $\kp_{(y,s),\ve}(x,t),\;\forall(x,t)\in \overline{\Om}_T,$ will stand for one of three functions $\al_{y,\ve},\;\beta_{y,\ve}$(if $s=0$) and $\gamma_{(y,s),\ve}$(if $s>0$), see (\ref{sub5}), (\ref{sub6}) and (\ref{sub11}). Thus, for every $(y,s)\in P_T$ and any $\ve>0$, small,
\eqRef{Exist1}
\kp_{(y,s),\ve}\in C(\overline{\Om}_T),\;\;\kp_{(y,s),\ve}(x,t)\le h(x,t),\;\forall(x,t)\in P_T,\;\;\mbox{and}\;\;\kp_{(y,s),\ve}(y,s)=h(y,s)-2\ve.
\ee
Moreover, $\kp_{(y,s),\ve}$ is a sub-solution in $\Om_T$, see (\ref{sub5}), (\ref{sub6}) and (\ref{sub11}) (see Part I of Appendix I for more details).
In an analogous manner, for every $(y,s)\in P_T$ and any $\ve>0$, small, $\kph_{(y,s),\ve}$ will denote one of the functions $\alh_{y,\ve},\;\bh_{y, \ve}$(if $s=0$) and 
$\gh_{(y,s),\ve}$(if $s>0$), see (\ref{sup2}), (\ref{sup3}) and (\ref{sup10}). Thus,
\eqRef{Exist2}
\kph_{(y,s),\ve}\in C(\overline{\Om}_T),\;\;\kph_{(y,s),\ve}(x,t)\ge h(x,t),\;\forall(x,t)\in P_T,\;\;\mbox{and}\;\;\kph_{(y,s),\ve}(y,s)=h(y,s)+2\ve.
\ee
Also, $\kph_{(y,s),\ve}$ is a super-solution in $\Om_T$ (see Part II of Appendix I for more details).
\vsp
\NI We adapt the ideas described in \cite{CI} and \cite{CIL} in what follows. Let $u,\;v:P_T\cup \Om_T\rightarrow \IR^+$. Suppose that $u$ is upper semi-continuous in $P_T\cup \Om_T$ and $v$ is lower semi-continuous in $P_T\cup \Om_T$. Assume further that $u(x,t)\le h(x,t)\le v(x,t),\;\forall(x,t)\in P_T$. 
Define
\eqRef{Exist3}
\aleph_{sub}=\{u(x,t):\mbox{$u$ is a sub-solution of (\ref{Exist0})} \}
\;\mbox{and}\;\aleph_{sup}=\{v(x,t):\mbox{$v$ is a super-solution of (\ref{Exist0})} \}.
\ee
See (\ref{def8}) and (\ref{def9}). We record that for every $(y,s)\in P_T$ and any $\ve>0$, small, $\kp_{(y,s),\ve}\in \aleph_{sub}$ and $\kph_{(y,s),\ve}\in \aleph_{sup}$.
Set
\eqRef{Exist31}
u_{sub}(x,t)=\sup_{u\in \aleph_{sub}} u(x,t), \;\;\mbox{and} \;\; v_{sup}(x,t)=\inf_{v\in \aleph_{sup}}  v(x,t),\;\;\forall(x,t)\in \Om_T.
\ee
We observe the following:

\begin{lem}\label{Exist4}
Let $\Om\subset \IR^n$ be a bounded domain and $T>0$. Suppose that $\aleph_{sub}$ and $\aleph_{sup}$ are as defined in (\ref{Exist3}). Set
$$u_{sub}(x,t)=\sup_{u\in \aleph_{sub}} u(x,t) \;\;\mbox{and} \;\; v_{sup}(x,t)=\inf_{v\in \aleph_{sup}}  v(x,t),\;\;\forall(x,t)\in \Om_T\cup P_T$$
Then (i) $u_{sub}\le v_{sup}$ in $\Om_T$, (ii) $u_{sub}(x,t)=v_{sup}(x,t)=h(x,t),\;\forall(x,t)\in P_T$, and for any $(y,s)\in P_T$,
$$(iii)\quad lim_{(x,t)\rightarrow (y,s)} u_{sub}(x,t)=\lim_{(x,t)\rightarrow (y,s)} v_{sup}(x,t)=h(y,s),\;\;\mbox{where}\;(x,t)\in \Om_T.$$
\end{lem}
\NI{\bf Proof:} If $u\in \aleph_{sub}$ and $v\in \aleph_{sup}$ then by Theorem \ref{comp}, $u\le v$ in $\Om_T$, since $u\le h\le v$ on $P_T$. Thus, $u_{sub}\le v$, in $\Om_T$ for any $v\in \aleph_{sup}$ and (i) holds. 
To show (ii) and (iii), let $(y,s)\in P_T$, and $\ve>0$, small. By (\ref{Exist1}), (\ref{Exist2}) and (\ref{Exist3}),
$$\kp_{(y,s),\ve}(x,t)\le u_{sub}(x,t)\le v_{sup}(x,t)\le \kph_{(y,s),\ve}(x,t),\;\forall(x,t)\in P_T\cup \Om_T.$$
Taking $(y,s)=(x,t)$, and since $\ve$ is arbitrary we obtain (ii). By letting $(x,t)\rightarrow (y,s)$, $(x,t)\in \Om_T$, in the above, we have
$$h(y,s)-2\ve\le \liminf_{(x,t)\rightarrow (y,s)}u_{sub}(x,t)\le \limsup_{(x,t)\rightarrow (y,s)} v_{sup}(x,t)\le h(y,s)+2\ve.$$
This yields (iii). $\Box$
\vsp
\NI Recall from (\ref{def3}) the definition of the cylinder $D_{\dl,\dl}(x,t)=B_\dl(x)\times(t-\dl/2,t+\dl/2).$ We now introduce the following definitions, see \cite{CIL}. 
For $b:\Om_T\cup P_T\rightarrow \IR$, we define the upper semi-continuous envelope of $b$ as
\eqRef{Exist41}
b^{us}(z,\tht)=\lim_{\dl\downarrow 0} \left( \sup\{b(x,t):\;(x,t)\in (\Om_T\cup P_T)\cap D_{\dl,\dl}(z,\tht)\} \right),\;\;\forall(z,\tht)\in \Om_T\cup P_T.
\ee
Similarly, we define the lower semi-continuous envelope of $b$ as
$$b^{\ell s}(z,\tht)=\lim_{\dl\downarrow 0} \left( \inf\{b(x,t):\;(x,t)\in (\Om_T\cup P_T)\cap D_{\dl,\dl}(z,\tht)\} \right),\;\;\forall(z,\tht)\in \Om_T\cup P_T.$$
We show that $u_{sub}^{us}$ is a sub-solution and $v_{sup}^{\ell s}$ is a super-solution of (\ref{Exist0}). The proof is an adaptation of Proposition 4.3 in \cite{CIL}, we provide details for completeness, also see \cite{CI}.

\begin{lem}\label{Exist5} Suppose that $w$ is upper semi-continuous on $\Om_T\cup P_T$, $(z,\tht)\in \Om_T$ and $(a,p,X)\in \mathcal{P}^+_{\Om_T}u(z,\tht)$. Assume that
for each $m=1,2,\cdots,$ there is a $w_m:\Om_T\cup P_T\rightarrow \IR$, $w_m$ upper semi-continuous, such that 

\NI (i) there exists a sequence $(z_{m},\tht_{m})\in \Om_T$ with $(z_{m},\tht_{m},w_{m}(z_{m},\tht_{m}))\rightarrow (z,\tht, w(z,\tht)),$\\
\NI (ii) if $(x,t)\in \Om_T$ and for $j=1,2,\cdots,$ $(x_j, t_j)\in \Om_T$ is such that $(x_j,t_j)\rightarrow (x,t)$ then $\limsup_{j\rightarrow \infty} w_j(x_j,t_j)\le w(x,t).$\\
\NI Then there is a subsequence $w_{m_k},\;k=1,2,\cdots$, and sequences $(y_{m_k}, s_{m_k})\in \Om_T$ and \\
$(a_{m_k},p_{m_k},X_{m_k})\in \mathcal{P}^+_{\Om_T} w_{m_k}(y_{m_k}, s_{m_k})$ such that
$$(y_{m_k},s_{m_k}, w_{m_k}(y_{m_k},s_{m_k}), a_{m_k}, p_{m_k}, X_{m_k})\rightarrow (z,\tht, w(z,\tht), a, p, X),\;\mbox{as $k\rightarrow \infty.$} $$
\end{lem}
\NI{\bf Proof:} Note that $\tht>0.$ Without any loss of generality, we take $z=o$. We also observe that $|a|\le \ve+ (a^2/\ve)$, for any $\ve>0$ and any $a$.
Let $\dl>0$ be small. Since $(a,p,X)\in \mathcal{P}^+_{\Om_T}w(o,\tht)$ , there is an $r>0$, small, such that $\overline{D}_{r,r}(o,\tht)\subset \Om_T$ and 
$$w(x,t)\le w(o,\tht)+a(t-\tht)+\langle p, x\rangle+\frac{\langle X x,x\rangle}{2}+\dl ( |x|^2+ |t-\tht| ),\;\;\forall (x,t)\in \overline{D}_{r,r}(o,\tht),$$
see (\ref{def4}).
Next, for any $\ve$ with $0<\ve<r^2/16$, and $\forall(x,t)\in \overline{D}_{r,r}(o,\tht)$,
\bea\label{Exist6}
w(x,t) \le w(o,\tht)+a(t-\tht)+\langle p, x\rangle+\frac{\langle X x,x\rangle}{2}+\dl \left( |x|^2+ \frac{ (t-\tht )^2}{\ve}+\ve \right).
\eea 
For each $\ell=1,2,\cdots,$ let $(y_\ell,s_\ell)\in \overline{D}_{r,r}(o,\tht)$ be a point of maximum of 
\eqRef{Exist61}
F_\ell(x,t)=w_{\ell}(x,t)-\left( w(o,\tht)+a(t-\tht)+\langle p, x\rangle+\frac{\langle X x,x\rangle}{2}+2\dl \left[ |x|^2+ \frac{|t-\tht|^2}{\ve}+\ve \right]\right).
\ee
Note that $(y_\ell,s_\ell)$ depends on $\ve$. Since $F_\ell(x,t)\le F_\ell(y_\ell,s_\ell)$, rearranging terms, we have, $\forall(x,t)\in \overline{D}_{r,r}(o,\tht)$, 
\bea\label{Exist7}
w_{\ell}(x,t) \le w_{\ell}(y_\ell,s_\ell)&+&a(t-s_\ell)+\langle p, x-y_\ell\rangle+\frac{\langle Xx,x\rangle}{2}-\frac{\langle Xy_\ell,y_\ell\rangle}{2}\\
&+&2\dl \left( |x|^2-|y_\ell |^2+ \frac{(t-\tht)^2}{\ve} - \frac{ (s_\ell-\tht)^2}{\ve} \right).   \nonumber
\eea
Since $(z_\ell,\tht_\ell)\rightarrow (o,\tht)$, $(z_\ell,\tht_\ell)\in \overline{D}_{r,r}(o,\tht)$, if $\ell$ is large enough. For large $\ell$, take $(x,t)=(z_{\ell},\tht_{\ell})$, in (\ref{Exist7}), to obtain
\bea\label{Exist8}
w_{\ell}(z_{\ell},\tht_{\ell}) \le  w_{\ell}(y_{\ell},s_{\ell}) &+&a(\tht_{\ell}-s_\ell) +\langle p, z_{\ell}-y_\ell\rangle
+ \frac{\langle Xz_{\ell},z_{\ell}\rangle}{2}- \frac{\langle X y_\ell,y_\ell\rangle}{2} \\
&+& 2\dl \left( |z_{\ell}|^2-|y_\ell|^2+\frac{(\tht_{\ell}-\tht )^2}{\ve}-\frac{(s_\ell-\tht )^2}{\ve} \right). \nonumber
\eea
Working with a subsequence if needed, $(y_\ell,s_\ell)\rightarrow (y,s)\in \overline{D}_{r,r}(o,\tht)$. Note that $(y,s)$ depends on $\ve$. Next, (\ref{Exist8}) and hypotheses (i)  and (ii) of the lemma imply that
\bea\label{Exist81}
w(o,\tht) &\le& \liminf_{\ell\rightarrow \infty}w_{\ell}(y_{\ell},s_{\ell}) +a(\tht-s)-\langle p, y\rangle-\frac{\langle Xy,y\rangle }{2}-2\dl \left( |y|^2+\frac{(s-\tht)^2}{\ve} \right) \nonumber\\
&\le & w(y,s)+a(\tht-s)-\langle p, y\rangle-\frac{\langle Xy,y\rangle }{2}-2\dl \left( |y|^2+\frac{(s-\tht)^2}{\ve} \right).
\eea
Taking $(x,t)=(y,s)$ in (\ref{Exist6}), 
\eqRef{Exist82}
w(o,\tht)\ge w(y,s)+a(\tht-s)-\langle p, y\rangle-\frac{\langle Xy,y\rangle }{2}-\dl \left( |y|^2+\frac{(s-\tht)^2}{\ve}+\ve \right).
\ee
Combining this with (\ref{Exist81}), we get $|y|^2+[(s-\tht)^2/\ve]\le \ve$, and thus, for any $\ve>0$,
\eqRef{Exist83}
|y|\le \sqrt{\ve}\;\;\;\mbox{and}\;\;\;|s-\tht|\le \ve.
\ee
Since $\ve<r^2/16$, $(y,s)\in \overline{D}_{r/2,r/2}(o,\tht)$.
\vsp
\NI We now use a diagonalization type argument. Let $k_0$ be a natural number such that $1/k_0<r^2/16.$ For each $k=k_0,k_0+1,\cdots$, take $\ve=1/k$ in (\ref{Exist6}). For $k\ge k_0$, there are $M_k\ge 1$ and $L_k\ge 1$ such that $|z_m|^2+|\tht_m-\tht|\le 2/k$ and $|y_\ell|^2+|s_\ell-\tht|\le 2/k$, for $m\ge M_k$ and $\ell\ge L_k.$ The first follows from the hypothesis (i) of the lemma and the second from 
(\ref{Exist83}). Note that the sequence $(y_\ell,s_\ell)$ depends on $k$.
Taking $\ell,\;m>\max(M_k,L_k)$, we select an element of the sequence $\{(z_m,\tht_m)\}$ and an element of the sequence $\{(y_\ell,s_\ell)\}$ with $m=\ell$. Call this index $m_k$. Note that both $(z_{m_k},\tht_{m_k})$ and $(y_{m_k},s_{m_k})$ are in $D_{r,r}(o,\tht)$. From (\ref{Exist7}), for each $k=1,2\cdots$ and 
$\forall(x,t)\in \overline{D}_{r,r}(o,\tht)$, we have
\bea\label{Exist84}
w_{m_k}(x,t) \le w_{m_k}(y_{m_k},s_{m_k})&+&a(t-s_{m_k})+\langle p, x-y_{m_k}\rangle+\frac{\langle Xx,x\rangle}{2}-\frac{\langle Xy_{m_k},y_{m_k}\rangle}{2}  \nonumber\\
&+&2\dl \left( |x|^2-|y_{m_k} |^2+k (t-\tht)^2 -k (s_{m_k}-\tht)^2 \right).
\eea
Taking $(x,t)=(z_{m_k},\tht_{m_k})$, we get
\ben
w_{m_k}(z_{m_k},\tht_{m_k}) &\le& w_{m_k}(y_{m_k},s_{m_k})+a(z_{m_k}-s_{m_k})+\langle p, z_{m_k}-y_{m_k}\rangle+\frac{\langle Xz_{m_k},z_{m_k}\rangle}\nonumber\\
&-&\frac{\langle Xy_{m_k},y_{m_k}\rangle}{2} 
+2\dl \left( |z_{m_k}|^2-|y_{m_k} |^2+k (\tht_{m_k}-\tht)^2 -k (s_{m_k}-\tht)^2 \right).
\een
Letting $k\rightarrow \infty$, using (\ref{Exist83}), i.e., $\max( |s_{m_k}-\tht|,\;|\tht_{m_k}-\tht|,\;|y_{m_k}|^2 )\le 2/k$ and the hypothesis (i) of the lemma, we have
$$w(o,\tht)\le \liminf_{k\rightarrow \infty} w_{m_k}(y_{m_k},s_{m_k})\le \limsup_{k\rightarrow \infty}w_{m_k}(y_{m_k},s_{m_k})\le w(o,\tht).$$
Next, we write 
\ben
(t-\tht)^2&=&(s_{m_k}-\tht)^2+2(t-s_{m_k})(s_{m_k}-\tht)+(t-s_{m_k})^2,\\
|x|^2&=&|(x-y_{m_k})+y_{m_k}|^2=|x-y_{m_k}|^2+2\langle x-y_{m_k}, y_{m_k}\rangle+|y_{m_k}|^2,\\
\langle Xx,x\rangle&=&\langle X(x-y_{m_k})+y_{m_k}, (x-y_{m_k})+y_{m_k}\rangle=\langle X(x-y_{m_k}), x-y_{m_k}\rangle\\
&+& 2\langle Xy_{m_k}, x-y_{m_k}\rangle+ \langle Xy_{m_k},y_{m_k}\rangle. 
\een
Using these in (\ref{Exist84}) we obtain
\ben 
w_{m_k}(x,t) \le  w_{m_k}(y_{m_k},s_{m_k})& + &(a+4\dl k(s_{m_k}-\tht))(t-s_{m_k}) + \langle p + 4\dl y_{m_k}+ X y_{m_k}, x-y_{m_k} \rangle\\
&+& \frac{ \langle (X+4\dl I )(x-y_{m_k}), x-y_{m_k} \rangle }{2}+2\dl k(t-s_{m_k})^2.
\een
Then $(a+4k\dl(s_{m_k}-\tht), p+4\dl y_{m_k}+X y_{m_k}, X+4\dl I)\in \mathcal{P}^+_{\Om_T} w_{m_k}(y_{m_k}, s_{m_k}).$ By our selection of $s_{m_k}$, $k|s_{m_k}-\tht|\le 2,$ see (\ref{Exist83}). Working with a further subsequence, if needed, 
$(a+4k\dl(s_{m_k}-\tht), p+4\dl y_{m_k}+X y_{m_k}, X+4\dl I)\rightarrow (a+4c\dl, p, X+4\dl)$, for some appropriate $c$ with $-2\le c\le 2.$
Since the set of $(b,q,Y)\in \IR\times \IR^n\times S(n)$ such that there is a sequence $(x_j, t_j)\in \Om_T$, $(b_j, q_j, Y_j)\in \mathcal{P}^+_{\Om_T}\hat{w}_j(x_j,t_j)$
such that $(x_j,t_j, \hat{w}_j(x_j,t_j), b_j, q_j,Y_j)\rightarrow (o,\tht, w(o,\tht),b, q, Y)$ is closed, it contains $(a+4c\dl ,p,X+4\dl I)$. $\Box$.

\begin{rem}\label{Exist9} An analogous version of Lemma \ref{Exist5} holds when $w$ is lower semi-continuous. We work with $-w$. 
$\Box$
\end{rem}

\NI Recall the definitions of $\aleph_{sub}$, $\aleph_{sup}$, $u_{sub}$ and $v_{sup}$ from (\ref{Exist3}), (\ref{Exist31}) and Lemma \ref{Exist4}. Let $u^{us}_{sub}$ be the upper semi-continuous envelope of $u_{sub}$, see (\ref{Exist41}). 

\begin{lem}\label{Exist10} Let $u^{us}_{sub}$ be defined by (\ref{Exist3}), (\ref{Exist31}) and (\ref{Exist41}). Then $u^{us}_{sub}$ is a sub-solution of (\ref{Exist0}) and
$u^{us}_{sub}(x,t)=h(x,t),\;\forall(x,t)\in P_T$. Similarly, $v^{ls}_{sup}$ is a super-solution of (\ref{Exist0}) and $v^{ls}_{sup}(x,t)=h(x,t),\;\forall(x,t)\in P_T$. If $(z,\tht)\in P_T$ then
$$\lim_{(x,t)\rightarrow (z,\tht)}u^{us}_{sub}(x,t)=\lim_{(x,t)\rightarrow (z,\tht)}v^{ls}_{sup}=h(z,\tht).$$
Moreover, $u^{us}_{sub}\le v^{ls}_{sup}$ in $\Om_T$.
\end{lem}
\NI{\bf Proof:} The proof is as in Lemma 4.2 in \cite{CIL}. Here we use Lemma \ref{Exist5}, an adaptation of Proposition 4.3 in \cite{CIL}. Suppose that $(y,s)\in \Om_T$ and 
$(a,p,X)\in \mathcal{P}^+_{\Om_T} u^{us}_{sub}(y,s).$ Clearly, there is a sequence $(y_k,s_k)\in \Om_T$ and $u_k\in \aleph_{sub}$ such that 
$((y_k,s_k), u_k(y_k,s_k))\rightarrow ( (y,s), u^{us}_{sub}(y,s)).$ By Lemma \ref{Exist5}, there is a subsequence $((y_\ell,s_\ell), u_{\ell}(y_\ell,s_\ell))$ and 
$(a_\ell, p_\ell, X_\ell)\in  \mathcal{P}^+_{\Om_T}u_{\ell}(y_{\ell}, s_{\ell})$ such that
$$\langle Xp, p\rangle -3a(u^{us}_{sub}(y,s))^2=\lim_{\ell\rightarrow\infty} \langle X_\ell p_\ell, p_\ell \rangle-3a_\ell (u_\ell(y_\ell,s_\ell))^2\ge 0.$$  
The proof for $v^{ls}_{sup}$ is similar.
The rest of the lemma follows from Lemma \ref{Exist4} and Theorem \ref{comp}.  $\Box$
\vsp
\NI Call
\eqRef{Exist11}
\underline{u}=u^{us}_{sub}\;\;\;\mbox{and}\;\;\;\overline{v}=v^{ls}_{sup}.
\ee
By Lemma \ref{Exist10}, $\underline{u}^{ls}\le \overline{v}^{us}$ in $\Om_T$ and $\underline{u}^{ls}(x,t)=\overline{v}^{us}(x,t)=h(x,t),\;\forall(x,t)\in P_T$. Define
\eqRef{Exist12}
\aleph_{sol}=\{\phih(x,t):\;\mbox{$\underline{u}\le \phih\le \overline{v}$ and $\phih$ is a sub-solution of (\ref{Exist0}) } \}. 
\ee
Our goal is to show that $\phi(x,t)=\sup_{\phih\in \aleph_{sol}} \phih(x,t)$ is a solution to (\ref{Exist0}) or (\ref{pb1}). Before, we prove this, we need the following lemma (see Lemma 4.4 in \cite{CIL}), also see the calculations pertaining to this in Appendix II. 

\begin{lem}\label{Exist13} Suppose that $w:\Om_T\rightarrow \IR^+$ satisfies
$$\Df w\ge 3 w^2 w_t,\;\;\;\mbox{in $\Om_T$}.$$
Assume that $w^{ls}$ is not a super-solution, i.e., $\Df w^{ls}\le 3 (w^{ls})^2 (w^{ls})_t$ does not hold in $\Om_T$. Then for any small $\ve>0$, 
there is sub-solution $w_\ve$, i.e., $\Df w_\ve\ge 3 (w_\ve)^2(w_\ve)_t$, in $\Om_T$, such that
$$w_\ve(x,t)\ge w(x,t),\;\;\sup_{\Om_T}( w_\ve-w)>0,\;\;\mbox{and}\;\;\om_\ve(x,t)=w(x,t),\;\;\forall(x,t)\in \Om_T\setminus D_{\ve,\ve}(z,\tht).$$  
\end{lem}
\NI{\bf Proof:} Since $w^{ls}$ is not a super-solution, there is some $(z,\tht)\in \Om_T$ and $(a,p,X)\in \mathcal{P}^-_{\Om_T}w^{ls}(z,\tht)$ such that
$\langle Xp,p\rangle>3 a (w^{ls}(z,\tht))^2$. Without any loss of generality we can take $z=o$. Let $\mu=\langle Xp,p\rangle-3 a (w^{ls}(o,\tht))^2$.  
For $\dl>0$ and $\nu>0$, define
\eqRef{Exist130}
\psi_{\dl, \nu}(x,t)=w^{ls}(o,\tht)+a(t-\tht)+\langle p, x\rangle+\frac{ \langle X x, x\rangle}{2}+\dl-\nu( |x|^2+|t-\tht| ).
\ee
If $\nu$ and $\dl=\dl_0<\mu/2$ are small enough (depending on $a,\;p,\;X$ and $w(o,\tht)$), then by (\ref{app3}) (see Appendix II)
\ben
\Df \psi_{\dl_0,\nu}(x,t)-3(\psi_{\dl_0,\nu}(x,t))^2(\psi_{\dl_0,\nu})_t(x,t)>0,\;\forall(x,t)\in D_{\rh,\rh}(o,\tht),
\een
where $\rh>0$ is small enough. Since $(a,p,X)\in \mathcal{P}^-_{\Om_T}w^{ls}(o,\tht)$, we have 
$$w(x,t)\ge w^{ls}(x,t)\ge w^{ls}(o,\tht)+a(t-\tht)+\langle p, x\rangle +\frac{\langle Xx,x \rangle}{2}+o(|t-\tht|+|x|^2).$$ 
Selecting, $0<r\le \rh$, $\dl=\min\{\dl_0,\;(r^2/32)\nu \},$ and using (\ref{Exist130}), $w(x,t)>\psi_{\dl,\nu}(x,t),\;\forall(x,t)\in D_{r,r}(o,\tht)\setminus D_{r/2,(r/2)}(o,\tht),$ if $r$ is small enough. However, by the definition of $w^{ls}$, there are points near $o$, where $\psi_{\dl,\nu}>w$. Thus, there is a sequence of points $(z_m,\tht_m)\rightarrow (o,\tht)$ such that
$$\lim_{m\rightarrow \infty}(\psi_{\dl,\nu}(z_m,\tht_m)-w(z_m,\tht_m))=w^{ls}(o,\tht)+\dl-w^{ls}(o,\tht)>0.$$
Define
\ben
w_r(x,t)=\left\{ \begin{array}{lcr} \max\{ w(x,t),\;\psi_{\dl,\nu}(x,t)\}, && \forall(x,t)\in D_{r,r}(o,\tht),\\ w(x,t), && \forall(x,t)\in \Om_T\setminus D_{r,r}(o,\tht). \end{array} \right.
\een
Take $\ve=r$. The function $w_\ve$ solves $\Df w_\ve\ge 3(w_\ve)^2 (w_\ve)_t$ in $\Om_T$, by Lemma \ref{Exist5}. $\Box$.
\vsp
\NI We are now ready to state the existence theorem for (\ref{pb1}).

\begin{thm}\label{Exist14} Let $\Om\subset \IR^n$ be a bounded domain and $T>0.$ Let $\underline{u}$ and $\overline{v}$ be as in 
\eqref{Exist11} and $\aleph_{sol}$ be as in \eqref{Exist12}. Then
$$\phi(x,t)=\sup_{\phih\in \aleph_{sol}} \phih(x,t)$$
is $C(\Om_T\cup P_T)$ and is a solution of \eqref{Exist0} and \eqref{pb1}.
\end{thm}
\NI{\bf Proof:} The details are as in Theorem 4.1 in \cite{CIL}. It is clear that 
\eqRef{Exist15}
\underline{u}^{ls}\le \phi^{ls}\le \phi\le \phi^{us}\le \overline{v}^{us},\;\;\mbox{in $\Om_T\cup P_T$}.
\ee
Thus, $\phi(x,t)=h(x,t),\;\forall(x,t)\in P_T$. By Lemmas \ref{Exist5} and (\ref{Exist10}), $\phi^{us}$ is a sub-solution of (\ref{Exist0}) and by Lemma \ref{Exist4} and Theorem \ref{comp},
$\phi^{us}\le \overline{v}.$ Thus, $\phi^{us}\in \aleph_{sol}$. By definition, $\phi^{us}=\phi$ and $\phi$ is a sub-solution.  
If $\phi^{ls}$ is not a super-solution at some $(z,\tht)\in \Om_T$, then take $\phi_\ve$ as defined in Lemma \ref{Exist13}. If $\ve>0$ is small enough,
$\phi_\ve=h$ on $P_T$ and $\underline{u}\le \phi_\ve$. Since $\phi_\ve$ is a sub-solution, by Theorem \ref{comp}, $\phi_\ve\le \overline{v}$ and $\phi_\ve\in \aleph_{sol}$. Since $\phi_\ve>\phi$ some where in $\Om_T$, this contradicts the definition of $\phi$ and hence $\phi^{ls}$ is a super-solution. By Theorem \ref{comp}, $\phi\le \phi^{ls}$. By (\ref{Exist5}),
$\phi=\phi^{us}=\phi^{ls}$. Thus, $\phi$ is continuous and solves (\ref{Exist0}) or (\ref{pb1}). $\Box$
\vsp

\section{Appendix I: Sub-solutions and super-solutions for the problem in (\ref{pb1}) }
\vsp
\NI In this Appendix, we construct sub-solutions and super-solutions to (\ref{pb1}). These examples of sub-solutions and super-solutions will be arbitrarily close to the boundary data, in a local sense and prove important for the Perron method. 
\vsp
\NI We recall (\ref{pb1}) and Lemma \ref{pde}. The problem is to find a solution $\phi\in C(\Om_T\cup P_T)$ such that
\eqRef{bpb}
\Df \phi=3\phi^2 \phi_t,\;\mbox{in $\Om_T$, $\phi(x,0)=f(x),\;\forall x\in \Om$ and $\phi(x,t)=g(x,t),\;\forall(x,t)\in \p\Om\times[0,T),$} 
\ee
where $f\in C(\overline{\Om})$ and $g\in C(\p\Om\times[0,T))$. Define
$$m=\min\left[ \inf_{\Om} f(x),\;\inf_{\p\Om\times [0,T)} g(x,t) \right]\;\;\mbox{and}\;\;M=\max\left[ \sup_{\Om} f(x),\;\sup_{\p\Om\times [0,T)} g(x,t) \right].$$
We are interested in positive solutions and also assume that $0<m<M$.
In what follows, we will have occasion to use the change of variables indicated in Lemma \ref{pde}, i.e.,
$\phi>0,$ solves 
$$\Df \phi\ge \;(\le) \; 3\phi^2 \phi_t,\;\;\;\mbox{if and only if}\;\;\;\Df \eta+|D\eta|^4 \ge \;(\le)\; 3\eta_t,$$
where $\eta=\log \phi.$ We also recall the following notation for the cylinders in (\ref{def3}),
$$D_{r,\tau}(x,t)=B_r(x)\times(t-\tau,t+\tau)\;\;\;\mbox{and}\;\;\;F_{r,\tau}(x,t)=B_r(x)\times(t,t+\tau).$$
Define $h:P_T\rightarrow \IR^+$ as follows:
$$h(x,0)=f(x),\;\forall x\in \overline{\Om},\;\;\mbox{and} \;\;h(x,t)=g(x,t),\;\forall(x,t)\in \p\Om\times[0,T).$$

\begin{rem}\label{app}
Let $a\in \IR$ and $p\in \IR^n,\;n\ge 2.$ Suppose that $z\in \Om$ and $\tht>0$. \\
\NI (i) If $a(t-\tht)+o(|t-\tht|)\le (\ge)0,$ as $t\rightarrow \tht$ then $a=0$.\\
\NI(ii) If $\langle p, x-z\rangle+o(|x-z|)\le (\ge)0,$ as $x\rightarrow z$, then $p=0$.\\
\NI These follow quite easily. $\Box$
\end{rem}

\NI We set $\lam_R$ to be the first eigenvalue of $\Df$ on $B_R(o)$. These notations will be followed in the rest of this section. Next, we make an elementary observation that will be used in what follows.

\vsp
\NI{\bf Part 1: Sub-solutions} 
\vsp
\NI We divide our work into two steps. First, we tackle the initial data and then consider the data along $\p\Om\times [0,T)$. Some of the basic ideas used in constructing
the sub-solutions will also prove useful in constructing super-solutions.
\vsp
\NI In what follows, $\ve>0$ is small enough such that $m-2\ve>0$.
\vsp
\NI{\bf Step 1: Initial data} 
\vsp
\NI Let $y\in \Om$ and $r=|x-y|$. Since $h(x,0)$ is continuous on $\overline{\Om}$, there is a $0<\dl< $dist$(y,\p\Om)$, such that 
\eqRef{sub0}
h(y,0)-\ve \le h(x,0) \le h(y,0)+\ve,\;\; \forall x\in \overline{B}_{\dl}(y).
\ee
We assume that $h(y,0)>m$, see below.
By Theorem \ref{rad1}, for every $0<\lam<\lam_\dl$, there is a unique $\psi_\lam\in C(\overline{B}_\dl(y)),\;\psi_\lam(x)=\psi_\lam(r),$ such that
\eqRef{sub1}
\Df \psi_\lam+\lam \psi_\lam^3=0,\;\mbox{$\psi_\lam>0$, in $B_\dl(y),$ and $\psi_\lam=m-2\ve$ on $\p B_{\dl}(y)$.}
\ee
Also, $\psi_\lam(x)=\psi_\lam(r),\;0\le r\le \dl,$ is decreasing in $r$. Also, $\psi_\lam(0)$ is continuous in $\lam$.
Since $\psi_\lam(y)=\sup_{B_\dl(y)}\psi_\lam$, by Lemma \ref{bnd} and Remark \ref{rad2}, there is a $0<\lam<\lam_\dl$ such that
\eqRef{sub2}
\psi_\lam(y)=h(y,0)-2\ve.
\ee
Note that $\psi_\lam\ge m-2\ve$. Define $\phi_y(x,t)=m,\;\forall(x,t)\in \Om_T$, if $h(y,0)=m$. Otherwise, set
$$\phi_y(x,t)=\left\{  \begin{array}{lcr} \psi_\lam(x) e^{-\lam t/3}, && \forall(x,t)\in \overline{B}_\dl(y)\times[0,T),\\ (m-2\ve) e^{-\lam t/3}, && \forall(x,t)\in [\Om\setminus \overline{B}_\dl(y)]\times[0,T),\end{array} \right.
$$
where $\psi_\lam$ is as in (\ref{sub1}) and (\ref{sub2}). Note that $\phi_y(x,t)\ge (m-2\ve)e^{-\lam t/3},\;\forall(x,t)\in \Om_T,$ and $0<\phi_y(x,t)\le h(x,t),\;\forall(x,t)\in P_T$. We show next that $\phi_y(x,t)$ is a sub-solution.
\vsp
\NI Using Lemma \ref{sol1} and direct differentiation, $\phi_y$ is a sub-solution in $(B_\dl(y) \times (0,T)) \cup  ([\Om\setminus \overline{B}_\dl(y)]\times (0,T)) $. We check that 
$\phi_y$ is a sub-solution on $\p B_\dl(y)\times (0,T)$. Suppose that $(z,s)\in \p B_\dl(y)\times (0,T)$ and $(a,p,X)\in \overline{\mathcal{P}}^+_{\Om_T}\phi_y(z,s)$, i.e., for $(x,t)\rightarrow (z,s)$, 
with $(x,t)\in \Om_T$,
\eqRef{sub3}
\phi_y(x,t)\le \phi_y(z,s)+a(t-s)+\langle p, x-z\rangle+\frac{ \langle X (x-z), x-z\rangle }{2}+o(|t-s|+|x-z|^2),
\ee
see (\ref{def4}). Since $\phi_y(z,t)=(m-2\ve)e^{-\lam t/3},\;\forall t\in(0,T)$, taking $x=z$ in (\ref{sub3}) leads to
$(m-2\ve) \left( e^{-\lam t/3}- e^{-\lam s/3} \right) \le a(t-s)+o(|t-s|),$ as $t\rightarrow s.$
Hence,
\eqRef{sub4}
a=-\frac{(m-2\ve)\lam e^{-\lam s/3} } {3}<0.
\ee
Since $\phi_y(x,t)\ge (m-2\ve)e^{-\lam t/3}>0,\;\forall(x,t)\in \Om_T$, taking $t=s$ in (\ref{sub3}) we obtain
$$\langle p, x-z\rangle +o(|x-z|)\ge \phi_y(x,s)- (m-2\ve)e^{-\lam s/3}\ge 0,\;\;x\rightarrow z. $$
By Remark \ref{app}, $p=\vec{0}.$
Using (\ref{sub4}),
$\langle Xp, p\rangle=0\ge 3 a\phi_y^2(y,s) .$
\vsp
\NI We summarize as follows. Recalling (\ref{sub0}), (\ref{sub1}) and the calculations done above, we see that for every $y\in \Om$ and $\ve>0$, small,
\ben
&&(i)\; \phi_y\in C(\overline{\Om}_T),\; (ii) \;\Df \phi_y\ge 3\phi_y^2(\phi_y)_t,\;\mbox{in $\Om_T$}, \;(iii)\; 0<\phi_y(x,t)\le h(x,t),\;\forall(x,t)\in P_T,\\
&&\mbox{and}\;(iv)\; \phi_y(y,0)=h(y,0)-2\ve.
\een
The construction in (\ref{sub1}) and (\ref{sub2}) and its analogue will be used repeatedly in this section. The above construction excludes points that are on $\p\Om\times\{0\}$. These will be considered in Step 2. For later reference, set for every $y\in \Om$,
\eqRef{sub5}
\al_{y,\ve}(x,t)=\phi_y(x,t),\;\;(x,t)\in \overline{\Om}_T.
\ee

\NI{\bf Step 2:} We now consider the set $\p\Om\times [0,T).$
\vsp
\NI{\bf Case (a) $\p\Om\times\{0\}$:} Fix $y\in \p\Om$; set $r=|x-y|$. Since $h(x,t)$ is continuous in $P_T$, there is a $0<\dl=\dl(y)$ and $0<\tau=\tau(y)$ such that
\eqRef{sub51}
h(y,0)-\ve \le h(x,t)\le h(y,0)+\ve,\;\;\forall(x,t)\in (\overline{B}_\dl(y)\times [0,\tau])\cap P_T.
\ee
Recall (\ref{sub1}) and (\ref{sub2}) from Step 1. We assume that $h(y,0)>m$. Let $0<\lam<\lam_\dl$ be such that $\psi_\lam\in C(\overline{B}_\dl(y)),\;\psi_\lam>0,$ 
radial, solves
$\Df \psi_\lam+\lam \psi_\lam^3=0,$ in $B_\dl(y)$. Moreover, $\psi_\lam(r)$ is decreasing in $r$, $\psi_\lam(0)=h(y,0)-2\ve$, and $\psi_\lam|_{B_\dl(y)}=m-2\ve.$
\vsp
\NI We choose $k\ge \lam$ such that 
$$k\ge \frac{3}{\tau} \log \left ( \frac{h(y,0)-2\ve}{m-2\ve} \right). $$ 
We observe the following:
\ben
&&(i)\;\;0<\psi_\lam(x) e^{-kt/3}\le \psi_\lam(y) e^{-kt/3},\;\forall (x,t)\in (\overline{B}_\dl(y)\cap \p\Om)\times[0,T),\\
&&(ii)\;\;0<\psi_\lam(x) e^{-kt/3}\le (h(y,0)-2\ve)e^{-kt}\le h(x,t),\;\forall(x,t)\in (\overline{B}_\dl(y)\cap \p\Om)\times[0,\tau],\\
&&(iii)\;\;0<\psi_\lam(x) e^{-kt/3}\le\psi_\lam(y) e^{-k\tau/3}\le m-2\ve\le h(x,t),\;\forall (x,t)\in(\overline{B}_\dl(y)\cap \p\Om)\times[\tau,T).
\een
Using (\ref{sub51}), (ii) and (iii) above, we have
$$\psi_\lam(x) e^{-kt/3}\le h(x,t),\;\forall(x,t)\in \left[ \overline{B}_\dl(y)\times[0,\;T] \right] \cap P_T.$$ 
Define $\phi_y(x,t)=m$ in $\Om_T$, if $h(y,0)=m$, otherwise,
$$\phi_y(x,t)=\left\{ \begin{array}{lcr}\psi_\lam(x) e^{-kt/3}, & \forall(x,t)\in \left[ \overline{B}_\dl(y)\times (0,T) \right]\cap \Om_T,\\
                                                                  (m-2\ve) e^{-kt/3}, & \forall(x,t)\in \Om_T\setminus (\overline{B}_\dl(y)\times (0,T)). \end{array} \right. $$
By our construction, $\phi_y(x,t)\ge (m-2\ve)e^{-kt/3},\;\forall(x,t)\in \Om_T$. Using Lemma \ref{sol1}, we see that
\ben
\Df \phi_y-3\phi_y^2 (\phi_y)_t=e^{-kt} \left( \Df \psi_\lam+k\psi_\lam^3 \right)=
e^{-kt} \psi_\lam^3(k-\lam)\ge 0, \;\;\mbox{in $(B_\dl(y)\times(0,T))\cap \Om_T$.}
\een
A differentiation shows that $\phi_y$ is a sub-solution in $\Om_T\setminus (\overline{B}_\dl(y)\times(0,T))$. The proof that it is a solution along $(\p B_\dl(y)\times (0,T))\cap \Om_T$ is similar to Part 1. To recap, for every
$y\in \p\Om$ and $\ve>0$, small, $\phi_y(x,t)\le h(x,t)$ in $P_T$, $\phi_y(x,t)$ is a sub-solution in $\Om_T$, $\phi_y\in C(\overline{\Om}_T)$ and $\phi(y,0)=h(y,0)-2\ve$.
For a fixed $y\in \p\Om$, and any $\ve>0$, small, set 
\eqRef{sub6}
\beta_{y,\ve}(x,t)=\phi_y(x,t),\;\;\;\forall(x,t)\in \overline{\Om}_T.
\ee
\vsp
\NI{\bf Case (b) $\p\Om\times(0,T)$:} Let $(y,s)\in \p\Om\times (0,T)$. There is a $0<\dl_0=\dl_0(y,s)$ and $0<\tau_0=\tau(y,s)$ such that
\eqRef{sub7}
h(y,s)-\ve\le h(x,t)\le h(y,s)+\ve,\;\;\;\forall(x,t)\in \overline{D}_{\dl_0,\tau_0}(y,s)\cap P_T.
\ee
We will construct a sub-solution in a region $R\subset D_{\dl_0,\tau_0}$ and then extend it to the rest of $\Om_T$ as a sub-solution. The region $R$ will be the union of two opposing cones in $\IR^{n+1}$ along $t$, with their bases joined at $t=s$. See Figure 1. 

\begin{figure}[htbp]
\includegraphics[scale=.4]{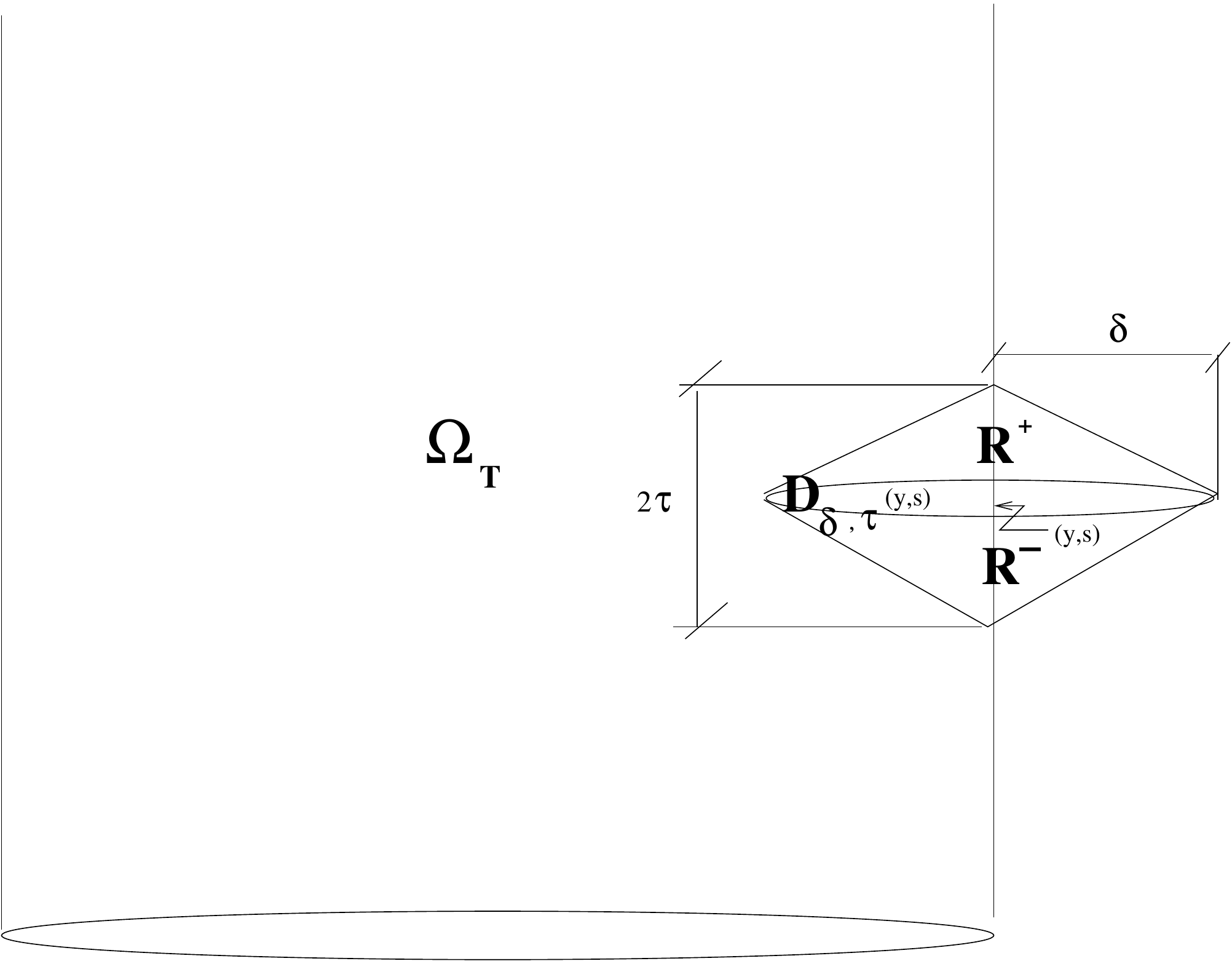}
\caption{Cones}
\label{fig:cones1}
\end{figure}

\NI In what follows, we select $0<\tau\le \tau_0$ and $0<\dl\le \dl_0$, constants $k$ and $c$ to ensure that we have a sub-solution. Also, to make our calculations easier, we will work with the inequality
$$\Df\eta+|D\eta|^4\ge 3\eta_t,\;\;\mbox{in $\Om_T$}$$
and $e^{\eta}$ will be the desired sub-solution for (\ref{pb1}), see Lemma \ref{pde}. Assume that $h(y,s)>m$.
\vsp
\NI We select
\eqRef{sub8}
k=\frac{1}{\tau} \log \left ( \frac{h(y,s)-2\ve}{m-2\ve} \right),\;\;\;c^4=3k,\;\;\;\mbox{and}\;\;\; \dl=\frac{k\tau}{c}.
\ee
A simple calculation gives us 
$$\dl=\left( \frac{\tau}{3} \right)^{1/4}  \left[ \log \left ( \frac{h(y,s)-2\ve}{m-2\ve} \right) \right]^{3/4}. $$
Since, $\dl\rightarrow 0$ as $\tau\rightarrow 0$, we can choose appropriate values of $\tau$ and $\dl$ such that $D_{\dl,\tau}(y,s)\subset D_{\dl_0,\tau_0}(y,s)$. This ensures that
(\ref{sub7}) is satisfied. The region $R$ will be in $D_{\dl,\tau}(y,s).$ 
We write $R=R^+\cup R^-$, where $R^+$ denotes the region in the upper cone and $R^-$ that in the lower cone. We describe these next. 
Set $r=|x-y|$, and choose $c,\;k,\;\tau$ and $\dl$ as in (\ref{sub8}). Define, in $0\le r\le \dl$ and $s-\tau\le t\le s+\tau$, 
\ben
&&R^+=\{(x,t):\;0\le cr\le k(s+\tau-t),\;s\le t\le s+\tau \},\;\;\mbox{and}\\
&& R^-=\{(x,t):\;0\le cr\le k(t-s+\tau),\;s-\tau\le t\le s \}.
\een
The apex of $R^+$ is $(y,s+\tau)$ and that of $R^-$ is $(y,s-\tau)$. Also, at $t=s$, $0\le r\le k\tau/c=\dl$. This is the set $R^+\cap R^-=\{(x,s):\;0\le |x-y|\le \dl\},$ the base common to both the cones.
\vsp
\NI If $h(y,s)=m$ then define $\etb(x,t)=m,\;\forall(x,t)\in \Om_T$, otherwise, first define the {\it{bump function} }
\eqRef{sub9}
\eta(x,t)=\left\{ \begin{array}{lcr} k(s+\tau-t)-cr+\log(m-2\ve), && \forall(x,t)\in R^+,\\ k(t-s+\tau)-cr+\log(m-2\ve), && \forall(x,t)\in R^- .
\end{array}\right.
\ee
The function $\eta$ is $C^{\infty}$ in the interior of $R^+\cap \Om_T$ and in $R^-\cap \Om_T$. Since $|x-y|>0,\;\forall x\in \Om$, differentiating and using (\ref{sub8}),
\ben
\Df \eta+|D\eta|^4-3\eta_t= \left\{ \begin{array}{lcr} c^4+3k\ge 0, && \forall(x,t)\in R^+\cap \Om_T,\\ c^4-3k=0, &&  \forall(x,t)\in R^-\cap \Om_T. \end{array} \right.
\een
\NI We now extend $\eta$ to the rest of $\overline{\Om}_T$, by
$$\bar{\eta}(x,t)=\left\{ \begin{array}{lcr} \eta(x,t), && \forall(x,t)\in R\cap \overline{\Om}_T,\\ \log(m-2\ve), && \forall(x,t)\in \overline{\Om}_T\setminus R.  \end{array} \right. $$
We make some observations. Recall that $R=R^+\cup R^-$. From (\ref{sub8}) and (\ref{sub9}), (i) $\eta(y,s)=\log(h(y,s)-2\ve)$, (ii) $\eta(x,t)=\log(m-2\ve),\;\forall(x,t)\in \p R$ and (iii) $\eta(x,t)\ge \log(m-2\ve), \forall(x,t)\in R$. 
Also, $0\le k(s+\tau-t)-cr\le k\tau$, in $R^+$, and $0\le k(t-s+\tau)-cr\le k\tau$, in $R^-$. Using (\ref{sub7}),
$$\log(m-2\ve)\le\eta(x,t)\le\log( h(y,s)-2\ve)\le \log(h(x,t)),\;\forall(x,t)\in R.$$
Hence, $\etb$ is continuous in $\overline{\Om}_T$, $\etb\ge \log(m-2\ve)$ in $\Om_t$ and $\etb\le h$ on $P_T$. Moreover, with $\dl$ as in (\ref{sub8}),
\eqRef{sub91}
\etb(x,s)=\log(m-2\ve), \forall x\in \p B_\dl(y)\cap \Om_T,\;\mbox{and}\;\etb(x,t)\le \log(h(y,s)-2\ve), \forall(x,t)\in R\cap P_T.
\ee
We show next that $\etb$ is a sub-solution in $\Om_T$.
To do this, we note that it is enough to check the definition on $\p R\cap \Om_T$ and the set common to both $R^+$ and $R^-$.  Recall that these lie in
the set $\overline{B}_\dl(y)\times[s-\tau,s+\tau]$. We describe these as follows.  \\
\NI (i) $\p R^+=\{(x,t)\in \IR^{n+1}:\;k(s+\tau-t)-c|x-y|=0,\;s< t\le s+\tau\}$,\\
\NI (ii) $\p R^-=\{(x,t)\in \R^{n+1}:\;k(t-s+\tau)-c|x-y|=0,\;s-\tau\le t< s\}$, and \\
\NI (iii)  $S=\{(x,s)\in\IR^{n+1}:\;0<|x-y|\le c\tau=\dl\}$. \\
The base $S$ common to both the cones is described in (iii).
\vsp
\NI We record the following. If $(z,\tht)\in (\p R^+\cup\p R^-\cup S)\cap \Om_T$ and $(a,p,X)\in \mathcal{P}^+_{\Om_T}\bar{\eta}(z,\tht)$ (recalling (\ref{def4})) then
for $(x,t)\rightarrow (z,\tht)$, $(x,t)\in \Om_T$,
$$\bar{\eta}(x,t) \le  \bar{\eta}(z,\tht) +a(t-\tht)+\langle p, x-z\rangle+\frac{\langle X(x-z),x-z\rangle}{2}+o(|t-\tht|+|x-z|^2).$$
Rearranging, we have, as $(x,t)\rightarrow (z,\tht)$, $(x,t)\in \Om_T$,
\eqRef{sub10}
a(t-\tht)+\langle p, x-z\rangle+\frac{\langle X(x-z),x-z\rangle}{2}+o(|t-\tht|+|x-z|^2)\ge \bar{\eta}(x,t)-\bar{\eta}(z,\tht).
\ee
In what follows, we will also make use of (\ref{sub91}) and the comment just preceding it. 
\vsp
\NI (i) Let $(z,\tht)\in\p R^+$, then $\tht>s$. By our construction, $\bar{\eta}(z,\tht)=\log(m-2\ve)$, $\bar{\eta}(x,t)\ge\bar{\eta}(z,\tht),\;\forall(x,t)\in \overline{\Om}_T$. If we take $x=z$ in (\ref{sub10}), then $a(t-\tht)+o(|t-\tht|)\ge 0,\;\mbox{as}\;\;t\rightarrow \tht.$ Using Remark \ref{app}, we get $a=0$. Take 
$t=\tht$ in (\ref{sub10}), and we get $\langle p, x-z\rangle+o(|x-z|)\ge 0$ as $x\rightarrow z$. We get $p=0$ from Remark \ref{app}. Thus, $\langle Xp,p\rangle+|p|^4=3a=0.$
\vsp
\NI (ii) Suppose that $(z,\tht)\in \p R^-$, then $\tht<s.$ Taking $x=z$ in (\ref{sub10}), using (\ref{sub8}) and (\ref{sub9}), we get
$a(t-\tht)+o(|t-\tht|)\ge\bar{\eta}(z,t)-\bar{\eta}(z,\tht)\ge 0 ,\;t\rightarrow \tht$. Remark \ref{app} implies that $a=0$. Next, taking $t=s$, we get
$\langle p, x-z\rangle+o(|x-z|)\ge 0,$ as $x\rightarrow z$. Thus, $p=0$ and $\langle Xp,p\rangle+|p|^4=3a=0.$
\vsp
\NI (iii) Now let $\tht=s$ and $(z,s)\in S\cap \Om_T$. If $(x,s)\in S\cap \Om_T$ then $x\in \overline{B}_\dl(y)$ and 
$$\bar{\eta}(x,s)=k\tau-c|x-y|+\log(m-2\ve).$$
 We consider the two sub-cases:
(1) $|z-y|=\dl$, and (2) $0<|z-y|<\dl.$ Call $\vec{e}_z=(z-y)/|z-y|$.
\vsp
\NI{\bf Sub-case (1):} Suppose that $|z-y|=\dl$, then $\etb(z,s)=\log(m-2\ve)$. Recall that $\etb(x,t)\ge \etb(z,s)$. In (\ref{sub10}), take $t=\tht=s$ to obtain $\langle p, x-z\rangle+o(|x-z|)\ge 0\;\mbox{as}\;\;x\rightarrow z,\;x\in \Om. $
Remark \ref{app} implies $p=0$. Next, take $x=z$ and notice that $a(t-\tht)+o(|t-\tht|)\ge \bar{\eta}(x,t)-\bar{\eta}(z,\tht)\ge 0,$ as $t\rightarrow 0$. This leads to $a=0$.
Thus, $\langle Xp,p\rangle+|p|^4=3a=0.$
\vsp
\NI{\bf Sub-case (2):} Set $\rh=|z-y|$ and let $0<\rh<\dl.$ Take $t=\tht=s$ in (\ref{sub10}) and use (\ref{sub8}) and (\ref{sub9}) to obtain
\eqRef{sub101}
c(\rh-|x-y|) \le \langle p, x-z\rangle+\frac{\langle X(x-z),x-z\rangle}{2}+o(|x-z|^2),\;\;x\rightarrow z,
\ee
since $\bar{\eta}(x,s)=k\tau-c|x-y|+\log(m-2\ve), \forall x\in B_\dl(y)$. Also, $\bar{\eta}$ is $C^2$ in $x$. Taking $|x-y|=\rho$, we get the spatial components of $p$ tangential to 
$\p B_{\rho}(y)$, at $z$, to be zero, and $p=\pm |p|\vec{e}_z$. If we take
$x=z+d \vec{e_z}$ with $d$, small, then $|x-y|=\rho+d$, and (\ref{sub101}) yields $-cd\le \pm |p|d+o(|d|)$ as $d\rightarrow 0$.
If $d<0$ we get $|p|\ge c$, $p=-|p|\vec{e_z}$, and $ -cd\le-|p|d+o(|d|)$ as $d\rightarrow 0$. Taking $d>0$, gives us $p=-c\vec{e}_z.$
Using the value of $p$, taking $x=z+d\vec{e}_z$ and recalling (\ref{sub101}), we have
$$0\le \frac{d^2}{2c^2} \langle X p,p\rangle +o(|d|^2),\;\;\mbox{as}\;\;x\rightarrow 0.$$
Thus,
$$\langle Xp, p\rangle\ge 0.$$
In (\ref{sub10}), take $\tht=s$ and $x=z$ to obtain $\bar{\eta}(z,t)-\bar{\eta}(z,s)\le a(t-s)+o(|t-s|)$, as $t\rightarrow s$. Using (\ref{sub9}),
$$a(t-s)+o(|t-s|)\ge \left\{ \begin{array}{lcr} k(t-s), && t\le s,\\  k(s-t), && t\ge s, \end{array} \right.\;\;\;\;\mbox{as}\;\;t\rightarrow s.$$
Hence, $-k\le a\le k,$ and recalling (\ref{sub8}) and that $|p|=c$, 
$$\langle Xp, p\rangle+ |p|^4-3a\ge c^4-3k=0.$$
Thus, $\bar{\eta}$ is a sub-solution and so $\phi=e^{\bar{\eta}}$ is a sub-solution to (\ref{bpb}). From (\ref{sub7}), (\ref{sub8}), (\ref{sub9}) and (\ref{sub91}), we have
$\phi(x,t)\le h(x,t),\;\forall(x,t)\in P_T$. To summarize, 
\ben
&&(i)\;\;\phi(x,t)\le h(x,t),\;\forall(x,t)\in P_T,\;(ii)\;\phi\in C(\overline{\Om}_T),\;(iii)\;\Df \phi\ge 3\phi^2\phi_t,\;\mbox{in $\Om_T$, and}\\
&&(iv)\; \phi(y,s)=h(y,s)-2\ve.
\een
Define for any $(y,s)\in \p\Om\times(0,T)$ and $\ve>0$,
\eqRef{sub11}
\gamma_{(y,s),\ve}=\phi(x,t)=e^{\eta(x,t)}. \;\;\;\Box
\ee
\vsp
\NI The three functions defined in (\ref{sub5}), (\ref{sub6}) and (\ref{sub11}) will be utilized to construct a sub-solution to (\ref{pb1}) and (\ref{Exist0}) that agrees with $h$ on $P_T$.
Let $u:P_T\cup\Om_T\rightarrow \IR^+$ be upper semi-continuous and $0<u(x,t)\le h(x,t),\;\forall(x,t)\in P_T$. We define
\eqRef{sub12}
\aleph_{sub}(\Om_T)=\{u(x,t):\;\mbox{$u$ is a sub-solution to (\ref{Exist0})} \}.
\ee
Then for every $\ve>0$, $\aleph_{sub}$ contains the functions $\al_{y,\ve},\;\beta_{y,\ve}$ and $\gamma_{(y,s),\ve}$, see (\ref{sub5}), (\ref{sub6}) and (\ref{sub11}). $\Box$

\begin{figure}[htbp]
\includegraphics[scale=.4]{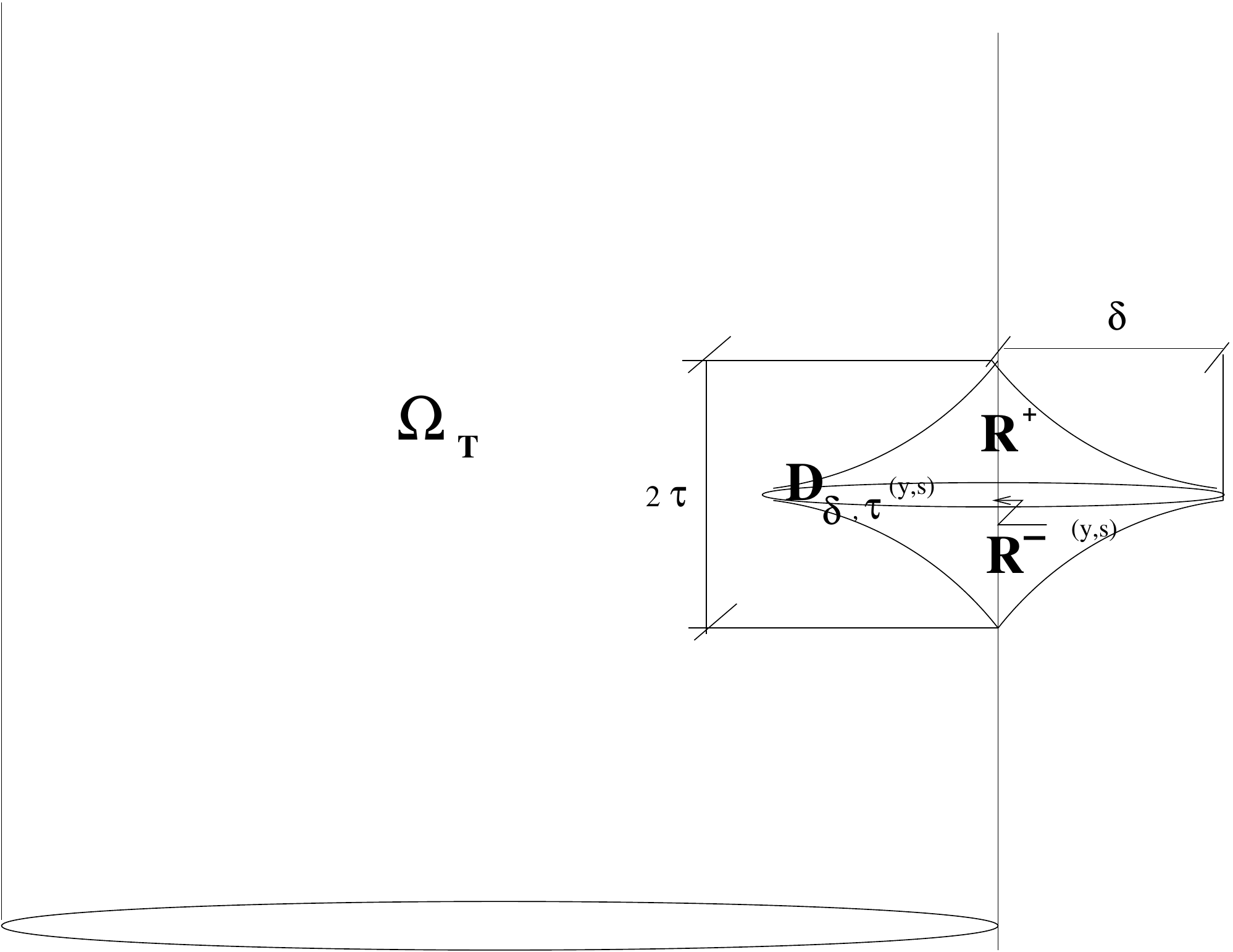}
\caption{Cusp shaped cones.}
\label{fig:cones2}
\end{figure}

\vsp
\NI{\bf Part 2: Super-solutions}
\vsp
\NI Next, we construct super-solutions to (\ref{bpb}). Our approach will be similar to what was done in Part 1.
Let $M=\sup_{(x,t)\in P_T} h(x,t)$. As done in Part 1, we divide our work into two steps. Let $\ve>0$ be small.
\vsp
\NI{\bf Step 1: Initial Data}
\vsp
\NI Let $y\in \Om$. Select $0<\dl<$dist$(y,\p\Om)$ such that $h(y,0)-\ve\le h(x,0)\le h(y,0)+\ve,\;\forall x\in \overline{B}_\dl(y)\cap \Om$. By Lemma \ref{rad4} and Remark \ref{radp}, for every $\lam>0$, there is a
function $\psi_\lam\in C(\overline{B}_\dl(y)),\;\psi_\lam>0,$ $\psi_\lam$ radial, such that 
\eqRef{sup1}
\Df \psi_\lam=\lam \psi_\lam^3,\;\;\mbox{in $B_\dl(y),$ and $\psi_\lam(y)=h(y,0)+2\ve.$}
\ee
We assume that $h(y,0)<M$. By Remark \ref{radp}, one can find a $\lam$ such that $\psi_\lam(x)=M+2\ve$, for $\forall x\in \p B_\dl(y).$ Define 
$$\phi(x,t)=\psi_\lam(x) e^{\lam t/3},\;\;\forall(x,t)\in \overline{B}_\dl(y)\times[0,T).$$
In case, $h(y,0)=M$, define $\phi(x,t)=M,\;\forall(x,t)\in \Om_T$, otherwise,
$$ \phih(x,t)=\left\{ \begin{array}{lcc} \psi_\lam(x) e^{\lam t/3}, && \forall(x,t)\in B_{\dl}(y)\times [0,T),\\ (M+2\ve)e^{\lam t/3}, && \forall(x,t)\in \Om_T\setminus (B_\dl(y)\times [0,T)). \end{array} \right.  $$
Note that $\psi_\lam(x)$, being radial, is also an increasing function of $r=|x-y|$. Thus, $\phih(x,t)\le (M+2\ve) e^{\lam t/3},\;\forall (x,t)\in \Om_T$.
From (\ref{sup1}) and Lemma \ref{sol1}, $\phih$ is a super-solution in $B_\dl(y)\times(0,T)$ and in its exterior. 
We need to check that $\phih$ is a super-solution on $\p B_\dl(y)\times(0,T).$ Let $(z,\tht)\in \p B_\dl(y)\times (0,T)$ and 
$(a,p,X)\in \mathcal{P}^-_{\Om_T}\phih(z,\tht)$. As $(x,t)\rightarrow (z,\tht)$, $(x,t)\in \Om_T$, 
$$\phih(x,t)\ge \phih(z,\tht)+a(t-\tht)+\langle p, x-z\rangle+\frac{\langle X(x-z), x-z\rangle}{2}+o(|t-\tht|+|x-z|^2). $$
Note that $\hat{\phi}(x,t)=(M+2\ve)e^{\lam t/3},\forall(x,t)\in \p B_\dl(y)\times(0,T).$
Since $\phih(z,\tht)=(M+2\ve)e^{\lam\tht/3}$, taking $x=z$, we get $a(t-\tht)+o(|t-\tht|)\le (M+2\ve)\left ( e^{\lam t/3}-e^{\lam \tht/3} \right),$ as $t\rightarrow \tht$. This yields $a=\lam(M+2\ve)e^{\lam \tht/3}/3$. If we take $t=\tht$, then
$$\langle p, x-z\rangle+o(|x-z|)\le \phih(x,\tht)-\phih(z,\tht)=\phih(x,\tht)-(M+2\ve)e^{\lam \tht/3}\le 0,\;\;\mbox{as}\;\;x\rightarrow z.$$
By Remark \ref{app}, $p=0$ and 
$\langle Xp,p\rangle\le 3\phih^2(z,\tht) a.$ Hence, $\phih$ is a super-solution in $\Om_T$. Also, recall that 
$h(x,0)\le h(y,0)+\ve\le \phih(x,0),\;\forall x\in \overline{B}_\dl(y)$. 
\vsp
\NI To summarize:
(i) $\phih\ge h,\;\mbox{in $P_T$},$ (ii) $\phih(y,0)=h(y,0)+2\ve,$ (iii) $\phih\in C(\overline{\Om}_T),\;\mbox{and},$ (iv) $\Df \phih\le 3\phih^2\phih_t,\;\mbox{in $\Om_T$.}$
For every $y\in \Om$ and $\ve>0$, set
\eqRef{sup2}
\alh_{y,\ve}(x,t)=\phih(x,t),\;\;\forall(x,t)\in \Om_T.
\ee
\vsp
\NI{\bf Step 2:} We consider $\p\Om\times [0,T)$. We consider the two cases, (a) $(y,s)\in \p\Om\times \{0\}$, and (b) $(y,s)\in \p \Om\times (0,T).$
\vsp
\NI{\bf Case (a) $\p\Om\times \{0\}$:} Let $y\in \p\Om$, then there are $\tau>0$ and $\dl>0$, small, such that
$$h(y,0)-\ve\le h(x,t)\le h(y,0)+\ve,\;\;\forall(x,t)\in (\overline{B}_\dl(y)\times [0,\tau] ) \cap P_T.$$
As done in Step 1(also see Case (a) of Step 2 in Part I), we select $\lam>0$ so that $\psi_\lam$ solves
$$\Df \psi_\lam=\lam \psi_\lam^3,\;\;\mbox{in $B_\dl(y)$, $\psi_\lam(y)=h(y,0)+2\ve,$ and $\psi_\lam(x)=M+2 \ve,\;\forall x\in \p B_\dl(y)$.}$$
Recall that $\psi_\lam(x)=\psi_\lam(r),$ is increasing in $r$, where $r=|x-y|$. Next, choose $k\ge \lam$ so that $(h(y,0)+2\ve)e^{k\tau/3}\ge (M+ 2 \ve).$ 
If $h(y,0)=M$, define $\phih=M$ in $\Om_T$, otherwise,
$$\phih(x,t)=\left\{ \begin{array}{lcc} \psi_\lam(x) e^{kt/3}, && \forall(x,t)\in (B_{\dl}(y)\times [0,T))\cap \overline{\Om}_T,\\ (M+2\ve)e^{kt/3}, && \forall(x,t)\in \Om_T\setminus (B_\dl(y)\times [0,T)). \end{array} \right. $$
By our choice of $\dl$ and $\tau$, $h(x,0)\le h(y,0)+2\ve \le \psi_\lam(x,0),\;\forall x\in \overline{B}_\dl(y)\cap \overline{\Om},$ and 
$$\phih(x,t)=\psi_\lam(x)e^{kt/3}\ge \left\{ \begin{array} {lcr} h(y,0)+2\ve\ge h(x,t), && \forall(x,t)\in (\overline{B}_\dl\times(0,\tau])\cap P_T,\\
(h(y,0)+2\ve)e^{k\tau/3}\ge M+2\ve,&& \forall(x,t)\in(\overline{B}_\dl(y)\times[\tau,T])\cap P_T.   \end{array} \right.$$
Also, $\phih(x,t)\le (M+2\ve)e^{kt/3},\;\forall(x,t)\in \Om_T$. 
\vsp
\NI Using Lemma \ref{sol1}, $\phih$  is a super-solution in $B_\dl(y)\times (0,T)$ and in $\Om_T\setminus(\overline{B}_\dl(y)\times (0,T)).$  Suppose that $(z,\tht)\in (\p B_\dl(y)\times (0,T))\cap \Om_T$, and $(a,p,X)\in \mathcal{P}^-_{\Om_T}\phih(z,\tht),$ i.e., as $(x,t)\rightarrow (z,\tht)$, $(x,t)\in \Om_T$,
$$\phih(x,t)\ge \phih(z,\tht)+a(t-\tht)+\langle p, x-z\rangle+\frac{\langle X(x-z), x-z\rangle}{2}+o(|t-\tht|+|x-z|^2),$$
Since $\phi(z,\tht)=(M+2\ve)e^{k\tht/3}$, taking $x=z$, one gets $a=k(M+2\ve)e^{k\tht/3}/3$. Taking $t=\tht$, one obtains
$\langle p, x-z\rangle+o(|x-z|)\le \phih(x,\tht)-\phih(z,\tht)\le 0.$ As done in Step 1, one can show that $p=0$, and
$$\langle Xp,p\rangle\le 3\phih^2 a.$$
We summarize: (i) $\phih\ge h,\;\mbox{in $P_T$},$ (ii) $\phih(y,0)=h(y,0)+2\ve,$ 
(iii) $\phih\in C(\overline{\Om}_T),\;\mbox{and},$ (iv) $\Df \phih\le 3\phih^2\phih_t,\;\mbox{in $\Om_T$.}$
For $y\in \p\Om$ and $\ve>0$, small, set
\eqRef{sup3} 
\bh_{y,\ve}(x,t)=\phih(x,t),\;\;\forall(x,t)\in \overline{\Om}_T.  \;\;\Box
\ee
\vsp
\NI{\bf Case (b) $\p\Om\times(0,T)$:} Let $(y,s)\in \p\Om\times(0,T)$, choose $\dl_1>0$ and $\tau_1>0$ such that 
\eqRef{sup40}
h(y,s)-\ve\le h(x,t)\le h(y,s)+\ve,\;\;\forall(x,t)\in (\overline{B}_{\dl_1}(y)\times(s-\tau_1,s+\tau_1))\cap P_T.
\ee
We will construct $\eta(x,t)$ such that $\Df \eta+|D\eta|^4\le 3\eta_t$, and $\phih=e^{\eta}$ will be the desired super-solution. 
Recall that $M=\sup_{P_T} h(x,t)$ and $m=\inf_{P_T}h(x,t).$ Assume that $h(y,s)<M$. Set
\eqRef{sup4}
\Gamma=\log\left( \frac{M+2\ve}{m+2\ve} \right),\;\;\nu=\frac{1}{1+2\Gamma},\;\;\;\mbox{and}\;\;\;k=\frac{1}{\tau} \log \left( \frac{M+2\ve}{h(y,s)+2\ve} \right).
\ee
Note that $k\tau\le\Gamma.$ We require that
\eqRef{sup5}
0<\dl\le \nu\left( \frac{k^2\tau^3 \Gamma}{3} \right)^{1/4},\;\;\mbox{and}\;\;c=\frac{k\tau}{\dl^\nu}.
\ee 
First select $0<\tau\le \tau_1$, next, use (\ref{sup5}) to choose $\dl\le \dl_1$, and finally, (\ref{sup4}) determines $c$. 
Let $r=|x-y|$. We construct our super-solution in a region $R$ that is the union of two cusp-shaped regions, centered at $(y,s)$. See Figure \ref{fig:cones2}.

\NI Set
\ben
&& R^+=\{(x,t):\; 0\le cr^\nu\le k(s+\tau-t),\;s\le t\le s+\tau,\;0\le r\le \dl\},\;\;\mbox{and}\\
&& R^-=\{(x,t):\; 0\le c r^\nu \le k(t-s+\tau),\;s-\tau\le t\le s,\;0\le r\le \dl\}.
\een
Call $R=R^+\cup R^-$, then $R\subset \overline{B}_\dl(y)\times[s-\tau,s+\tau]$. Define
\eqRef{sup6}
\eta(x,t)= \left\{ \begin{array}{lcr}  cr^\nu-k(s+\tau-t)+\log(M+2\ve), && \forall(x,t)\in R^+,\\ cr^\nu-k(t-s+\tau)+\log(M+2\ve), && \forall(x,t)\in R^-.  \end{array} \right.
\ee
Then, $\eta(y,s)=\log(h(y,s)+2\ve)\le \eta(x,t),\;\forall(x,t)\in R$. Hence, by (\ref{sup40}), $\log h(x,t)\le \eta(x,t)\le\log( M+2\ve),\;\forall(x,t)\in R\cap P_T$. We show that $\eta$ is a super-solution in the interior of $R^+\cap \Om_T$ and in the interior of
$R^-\cap \Om_T$. We start with $R^-\cap \Om_T$ and differentiate to find
$\eta_r=c\nu  r^{\nu-1},\;\;\eta_{rr}=c\nu (\nu-1)  r^{\nu-2},$ and $\eta_t=-k.$ Thus,
\ben
\Df \eta+|D\eta|^4-3\eta_t&=&\nu^3 c^3 (\nu-1) r^{3\nu-4}+\nu^4 c^4 r^{4\nu-4}+3k.\\
&=&3k+c^3 \nu^4 r^{3\nu-4} \left( c r^\nu-\frac{1-\nu}{\nu} \right).
\een
Note that $0<r< \dl$, $|t-s|<\tau$ and $0<\nu<1.$
We use (\ref{sup4}), (\ref{sup5}) and $c\dl^\nu=k\tau\le \Gamma$ to obtain
\ben
\Df \eta+|D\eta|^4-3\eta_t &=&3k+c^3 \nu^4 r^{3\nu-4} \left( c r^\nu-2\Gamma \right) \le 3k+c^3 \nu^4 r^{3\nu-4} \left( c \dl^\nu-2\Gamma \right)\\
&\le & 3k-c^3 \nu^4 \left( \frac{\Gamma}{r^{4-3\nu}} \right) \le 3k-c^3\nu^4 \left( \frac{\Gamma}{\dl^{4-3\nu} } \right),\\
&=&3k- \nu^4\frac{(c\dl^\nu)^3 \Gamma}{\dl^4}  =3k-\nu^4\frac{k^3 \tau^3 \Gamma}{\dl^4}\le 0.
\een
Differentiating in the interior of $R^+\cap \Om_T$, we have
$$\Df \eta+|D\eta|^4-3\eta_t \le -3k-\nu^4\frac{k^3 \tau^3 \Gamma}{\dl^4}\le 0.$$
\NI We now define 
\ben
\etb(x,t)=\left\{  \begin{array}{lcr} \eta(x,t), && \forall(x,t)\in R\cap \Om_T,\\  \log(M+2\ve), && \forall(x,t)\in \Om_T\setminus R.  \end{array}  \right.
\een
Then $\etb$ is continuous in $\overline{\Om}_T$ and $\etb\le \log(M+2\ve)$ in $\Om_T$. Define
\ben
&& \p R^+=\{(x,t):\;cr^\nu=k(s+\tau-t),\;s< t\le s+\tau,\;r\ge 0\},\\
&& \p R^-=\{(x,t):\;cr^\nu=k(t-s+\tau),\;s-\tau\le t< s,\;r\ge 0\},\\
&& S=\{(x,t):\;t=s,\; cr^\nu\le k\tau,\;0\le r\le \dl\}.
\een
\NI In order to check that $\etb$ is a super-solution, we need check the definition only on the boundaries $\p R^+\cap \Om_T$, $\p R^-\cap \Om_T$ and on the interface $S\cap \Om_T.$ Let $(z,\tht)$ be in one of these sets and $(a,p,X)\in  \mathcal{P}^-_{\Om_T}{\etb}(z,\tht)$, i.e., as $(x,t)\rightarrow (z,\tht)$, $(x,t)\in \Om_T$,
\eqRef{sup7}
a(t-\tht)+\langle p, x-z\rangle+\frac{ \langle X(x-z), x-z\rangle}{2}+o(|t-\tht|+|x-y|^2) \le \etb(x,t)-\etb(z,\tht).
\ee
\vsp
\NI In parts (i) and (ii) below, if $(z,\tht)\in \p R^+ \cup \p R^-$ then using (\ref{sup6}) we get $\etb(z,\tht)=\log(M+2\ve)$ and $\etb(x,t)\le \etb(z,\tht),\;\forall (x,t)\in \Om_T.$
\vsp
\NI (i) Let $(z,\tht)\in\p R^+\cap \Om_T$. 
Taking $x=z$ in (\ref{sup7}) and using (\ref{sup6}), we get $a(t-\tht)+o(|t-\tht|)\le 0$, as $t\rightarrow \tht. $   
Remark \ref{app} implies $a=0$. Taking $t=\tht$ in (\ref{sup7}), we get $\langle p, x-z\rangle+o(|x-z|)\le 0$, as $x\rightarrow z$. 
This gives us $p=0$ and $\langle Xp, p\rangle+|p|^4-3a=0$.
\vsp
\NI (ii) Let $(z,\tht)\in \p R^-\cap \Om_T$. Select $x=z$ in (\ref{sup7}), then (\ref{sup6}) implies $a(t-\tht)+o(|t-\tht|)\le 0,$ as $t\rightarrow \tht$.   
Thus, $a=0$. Next, take $t=\tht$, to find $\langle p, x-z\rangle+o(|x-z|)\le 0$, as $x\rightarrow z$. Remark \ref{app} implies
$p=0$ and $\langle Xp, p\rangle+|p|^4-3a=0$.
\vsp
\NI (iii) Let $(z,\tht)\in S\cap \Om_T$. Here $\tht=s$, and there are two sub-cases to consider: (1) $|z-y|=\dl$, and (2) $0<|z-y|<\dl.$ We rewrite (\ref{sup7}) as
\eqRef{sup8}
a(t-s)+\langle p, x-z\rangle+\frac{ \langle X(x-z), x-z\rangle}{2}+o(|t-s|+|x-z|^2)\le \etb(x,t)-\etb(z,s)
\ee
as $(x,t)\rightarrow (z,s)$, $(x,t)\in \Om_T$.
\vsp
\NI{\bf Sub-case 1:} Suppose that $|z-y|=\dl$ then $\etb(z,s)=\log(M+2\ve)$ and $\etb(x,t)\le \etb(z,s),\;\forall(x,t)\in \Om_T.$ Taking $x=z$ in (\ref{sup8}) and using (\ref{sup6}), we get
$a(t-s)+o(|t-s|)\le 0$, for $t\rightarrow s$. This leads to $a=0$. Now take $t=s$ in (\ref{sup8}) to get
$$\langle p, x-z\rangle+o(|x-z|)\le \etb(x,s)-\etb(z, s)\le 0,\;\;x\rightarrow z.$$
Remark \ref{app} yields $p=0$. Thus, $\langle Xp, p\rangle+|p|^4-3a=0$.
\vsp
\NI{\bf Sub-case 2:} Suppose that $0<|z-y|<\dl.$ Note that $\etb\in C^2$ in $x$ in $B_\dl(y)$. For $\xi\in \IR^n$, call $d_\xi=|\xi-y|$. Set $t=s$ in (\ref{sup8}) and use (\ref{sup6}) to get
\eqRef{sup81}
\langle p, x-z\rangle +o(|x-z|)\le c(d_x^\nu-d_z^\nu).
\ee
Call $\vec{e_z}=(z-y)/d_z$. If $d_x=d_z$, the components of $p$ tangential to $\p B_{d_z}(y)$, at $z$, are zero. Thus, $p=\pm |p|\vec{e_z}$. Take $x=y+d_x \vec{e_z}$, then
$x-z=(d_x-d_z)\vec{e}_z$, and (\ref{sup81}) implies 
$\pm |p| ( d_x-d_z)+o(|d_x-d_z|)\le c(d_x^\nu-d_z^\nu),\;x\rightarrow z.$ This gives $p=|p|\vec{e}_z=c\nu d_z^{\nu-1}\vec{e_z}$. Using the second order expansion in (\ref{sup8}) ($t=s$),
$$c\nu d_z^{\nu-1}(d_x-d_z)+\frac{(d_x-d_z)^2}{2|p|^2} \langle Xp, p\rangle+o(|x-z|^2)\le c(d_x^\nu-d_z^\nu),\;\;x\rightarrow z.$$
Noting that $d_z>0$, a second order expansion on the right hand side of the inequality leads to
\eqRef{sup9}
\langle Xp, p \rangle\le c\nu (\nu-1)d_z^{\nu-2}|p|^2.
\ee
Next, we take $x=z$ in (\ref{sup8}) and use (\ref{sup6}) to see that
$$a(t-s)+o(|t-s|)\le \left\{ \begin{array}{lcr} k(t-s), && t\ge s,\\ k(s-t), && t\le s. \end{array} \right.$$
Thus, $-k\le a\le k$. Using (\ref{sup4}) (note $0<\nu<1$), (\ref{sup5}), (\ref{sup9}), $|p|= c\nu d_z^{\nu-1}$ and $cd_z^\nu\le c\dl^\nu\le \Gamma$, we calculate
\ben
\langle Xp, p\rangle+|p|^4-3a &\le& c\nu(\nu-1)d_z^{\nu-2}|p|^2+|p|^4+3k
=c^3\nu^3(\nu-1)d_z^{3\nu-4}+c^4\nu^4 d_z^{4\nu-4}+3k\\
&=&c^3 \nu^4 d_z^{3\nu-4}\left( cd_z^{\nu}-\frac{1-\nu}{\nu} \right)+3k
\le c^3 \nu^4 d_z^{3\nu-4}\left( c\dl^{\nu}-\frac{1-\nu}{\nu} \right)+3k\\
&\le & 3k-c^3 \nu^4 d_z^{3\nu-4}\Gamma=3k-c^3 \nu^4 \frac{\Gamma}{d_z^{4-3\nu}}\le 3k-c^3 \nu^4 \frac{\Gamma}{\dl^{4-3\nu}}\\
&\le& 3k-c^3 \nu^4 \dl^{3\nu}\frac{\Gamma}{\dl^{4}}=3k-\nu^4 (k\tau)^3\frac{\Gamma}{\dl^{4}}\le 0.
\een
Thus, $\etb$ is super-solution. We summarize (see (\ref{sup6}) and the definition of $\etb$):
\ben
&&(i)\;\etb\in C(\overline{\Om}_T),\;(ii)\; \etb (y,s)=h(y,s)+2\ve,\;(iii)\;\etb(x,t)\ge h(x,t),\;\forall(x,t)\in P_T, \;\mbox{and}\\
&&(iv)\; \Df \etb+|D\etb|^4\le 3\etb_t,\;\mbox{in $\Om_T$.}
\een
For every $(y,s)\in \p\Om\times(0,T)$ and $\ve>0$, small, define
\eqRef{sup10}
\gh_{(y,s),\ve}=e^{\eta},\;\;\mbox{in $\overline{\Om}_T$}.\;\;\;\Box
\ee
\vsp
\NI The three functions defined in (\ref{sup2}), (\ref{sup3}) and (\ref{sup10}) will be utilized to construct a super-solution to (\ref{pb1}) and (\ref{Exist0}) that agrees with $h(x,t)$ on $P_T$.
\NI Let $v:P_T\cup\Om_T\rightarrow \IR$ be lower semi-continuous and $v(x,t)\ge h(x,t),\;\forall(x,t)\in P_T$. We define
\eqRef{sup11}
\aleph_{sup}(\Om_T)=\{u(x,t):\;\mbox{$v$ is a super-solution in $\Om_T$} \}.
\ee
Then for every $\ve>0$, $\aleph_{sup}$ contains the functions $\alh_{y,\ve},\;\bh_{y,\ve}$ and $\gh_{(y,s),\ve}$, see (\ref{sup2}), (\ref{sup3}) and (\ref{sup10}).

\section{Appendix II}
\vsp
\NI {\bf Calculations for Lemma \ref{Exist13}. } Suppose that $k\in \IR$ and $(a,p,X)\in \IR\times \IR^n\times S(n)$ are such that $\langle Xp, p\rangle-3 a k^2>0.$ Let $\tht>0$, $\dl>0$ and $\nu>0$; define for $(x,t)$ in an open set, in $\IR^{n+1}$, containing $(o,\tht)$,
\eqRef{app0}
u(x,t)=k+a(t-\tht)+\langle p, x\rangle+\frac{ \langle X x, x\rangle}{2}+\dl-\nu( |x|^2+|t-\tht| ).
\ee
We claim that if $\dl$ and $\nu$ are small enough then $\Df u-3u^2u_t\ge 0$ in a small neighborhood of $(o,\tht)$. Let $(b,q,Y)\in \mathcal{P}^+u(y,s)$, i.e, as 
$(x,t)\rightarrow (y,s)$,
\eqRef{app1}
u(x,t)-u(y,s)\le b(t-s)+\langle q, x-y\rangle+\frac{ \langle Y(x-y),x-y\rangle }{2}+o(|t-s|+|x-y|^2).
\ee
Using (\ref{app0}) together with the inequality $|\; |t-\tht|-|s-\tht|\; |\le |t-s|$, we get
\bea\label{app2}
u(x,t)-u(y,s)&=&a(t-s)+\langle p, x-y\rangle+\frac{ \langle X x, x\rangle}{2}-\frac{ \langle X y, y\rangle}{2} \nonumber\\
&-&\nu( |x|^2-|y|^2+|t-\tht|-|s-\tht| )  \nonumber\\
&\ge & a(t-s)-\nu |t-s|+\langle p+Xy-2\nu y, x-y\rangle+ \frac{ \langle (X-2\nu I)(x-y), x-y\rangle}{2}.
\eea
Take $x=y$ in (\ref{app1}) and (\ref{app2}) to find $b(t-s)+o(|t-s|)\ge a(t-s)-\nu |t-s|$ as $t\rightarrow s$. This yields $a-\nu\le b\le a+\nu$. Next, take $t=s$ in (\ref{app1}) and (\ref{app2})
to obtain $q=p+Xy-2\nu y.$ Using this fact we see that $\langle Y(x-y),x-y\rangle +o(|x-y|^2)\ge \langle (X-2\nu I)(x-y),x-y\rangle$, for every $(x,t)\rightarrow (y,s)$. Thus, $Y\ge X-2\nu I$. Next, noting that $u(o,\tht)=k$, using (\ref{app0}) and $b\le a+\nu$,
\bea\label{app3}
&& \langle Y q,q\rangle -3 b u(y,s)^2 \ge \langle (X-2\nu I) q,q\rangle -3 b k^2+3b(k^2-u(y,s)^2)\\
&&=\langle X q,q\rangle -3 b k^2+O(|\nu|+|y|+|t-\tht|)\ge\langle Xq, q\rangle-3ak^2+O(|\nu|+|y|+|t-\tht|)\nonumber\\
&&= \langle X(p+Xy-2\nu y), (p+Xy-2\nu y)\rangle-3ak^2+O(|\nu|+|y|+|t-\tht|)\nonumber\\
&&=\langle Xp, p\rangle-3ak^2+O(|y|+|\nu|+|t-\tht|).\nonumber
\eea
If $|y|,\;|t-\tht|$ and $\nu$ are small enough, $u$ satisfies $\Df u-3 u^2u_t\ge 0,$ near $(o,\tht)$.

\vsp
\NI Department of Mathematics\hfil
Western Kentucky University\hfil
Bowling Green, Ky 42101
\vsp
\NI Department of Mathematics\hfil
University of Kentucky\hfil
Lexington, KY 40506


\begin{thebibliography}{99}

\bibitem{AS} G. Akagi and K. Suzuki, \it Existence and uniqueness of viscosity solutions for a degenerate parabolic equation associated with the infinity-Laplacian, \rm
Calc. Var. (2008), 31, 457-471. 

\bibitem{AJK} G. Akagi, P. Juutinen and R. Kajikiya, \it Asymptotic behavior of viscosity solutions for a degenerate parabolic equation associated with the infinity-Laplacian, \rm
Math. Annalen, 343 (2009), no 4, 921-953.

\bibitem{ACJ} G. Aronson, M. Crandall, P. Juutinen, \it A tour of the theory of absolute minimizing functions. \rm Bull. Amer. Math. Soc., 41 (2004), 439-505.

\bibitem{BDM} T. Bhattacharya, E. DiBenedetto and J. J. Manfredi, \it Limits as $p\rightarrow \infty$ of $\Delta_pu_p=f$ and related extremal problems, \rm Some topics in nonlinear PDEs(Turin, 1989), Rend. Sem. Mat.Univ. Politec. Torino 1989, Special Issue, 15-68 (1991).

\bibitem{BL} Tilak Bhattacharya and Leonardo Marazzi, \it An eigenvalue problem for the infinity-Laplacian, \rm Electronic Journal of Differential Equations, Vol. 2013(2013), no 47, 1-30.

\bibitem{BMO1} T. Bhattacharya and A. Mohammed, \it On solutions to Dirichlet problems involving the infinity-Laplacian, \rm Advances in Calculus of Variations, vol 4, issue 4, 445-487 (2011).            

\bibitem{BMO2} T. Bhattacharya and A. Mohammed, \it Inhomogeneous Dirichlet problems  involving the infinity-Laplacian, \rm Advances in Differential Equations, vol 17, nos 3-4,
225-266. 

\bibitem{C} M. G. Crandall, \it A visit with the $\infty$-Laplace equation, \rm  Calculus of Variations and Nonlinear Partial Differential Equations, Lecture Notes in Math. 1927.
pp 75-12, Springer Berlin, 2008.

\bibitem{CEG} M. G. Crandall, L. C. Evans and R. F. Gariepy, Optimal Lipschitz extensions and the infinity-Laplacian, Calc. Var. Partial Differential Equations 13(2001), no 2, 123-139.

\bibitem{CI} M. G. Crandall and H. Ishii, The Maximum Principle for Semicontinuous Functions, Differential and Integral Equations, Vol 3, No 6, 1990, 1001-1014.

\bibitem{CIL} M. G. Crandall, H. Ishii and P. L. Lions, \it User's guide to viscosity solutions of second order partial differential equations,\rm Bull. Amer. Math. Soc. 27(1992) 1-67.

\bibitem{CW} M. G. Crandall and P. Wang, \it Another way to say caloric, \rm Dedicated to Philippe Benilan. J. Evol. Equ 3 (2003), no 4, 653-672. 

\bibitem{ED} E. DiBenedetto, Degenerate Parabolic Equations, Universitext, Springer, 1993.

\bibitem{LS} L. C. Evans and C. K. Smart, \it Everywhere differentiability of infinity-harmonic functions. \rm To appear.

\bibitem{LSA} L. C. Evans and O. Savin, \it $C^{1,\al}$ regularity for infinity-harmonic functions in two variables,\rm Calc. Var. Partial Differential Equations, 32 (2008), 325-347.

\bibitem{LN} E. M. Landis, \it Second Order Equations of Elliptic and Parabolic Type,\rm Translations of Mathematical Monographs (AMS) 1998.

\bibitem{PSSW} Y. Peres, O. Schramm, S. Sheffield and D. Wilson, \it Tug of war and the infinity-Laplacian, \rm J. Amer. Math. Soc. 22 (2009), no 1, 167-210. 

\end{thebibliography}
\end{document}